\documentclass[11pt]{amsart}

\usepackage{amsfonts,latexsym,graphicx,fullpage}

\newtheorem{proposition}{Proposition}
\newtheorem{theorem}[proposition]{Theorem}
\newtheorem{lemma}[proposition]{Lemma}
\newtheorem{corollary}[proposition]{Corollary}
\newtheorem{conjecture}[proposition]{Conjecture}

\theoremstyle{definition}
\newtheorem{definition}[proposition]{Definition}

\theoremstyle{remark}

\newtheorem{example}[proposition]{Example}

\def\R{\mathbb{R}}
\def\Z{\mathbb{Z}}
\def\d{\partial}
\def\a{\tilde{a}}
\def\sgn{\mbox{sgn }}

\def\im{{\rm im}}
\def\A{\mathbb{A}}
\def\CC{\mathcal{C}}

\def\alg#1{\Z[t,t^{-1}]\langle{#1}\rangle}
\def\eps{\varepsilon}

\def\pixy{\pi_{xy}}

\numberwithin{proposition}{section}

\theoremstyle{plain}

\begin{document}

\title{Computable Legendrian Invariants}
\author{Lenhard L.\ Ng}
\address{American Institute of Mathematics, 360 Portage Avenue,
Palo Alto, CA 94306}
\curraddr{Department of Mathematics, Massachusetts Institute of
Technology, Cambridge, MA 02139}
\email{lenny@math.mit.edu \newline \indent
\textit{URL}: http://www-math.mit.edu/\~{}lenny/}
\thanks{The author was partially supported by the American
Institute of Mathematics Five-Year Fellowship, and by research
assistantships provided by grants from the NSF and DOE}
\renewcommand{\subjclassname}{\textup{2000} Mathematics Subject Classification}

%\renewcommand{\datename}{}
%\date{Comments are welcome; please send correspondence to
%\texttt{lenny@math.mit.edu}.}
%\address{\newline \indent Department of Mathematics \newline
%\indent Massachusetts Institute of Technology \newline
%\indent 77 Massachusetts Avenue \newline
%\indent Cambridge, MA 02139 \newline}
%\email{lenny@math.mit.edu \newline 
%\indent {\it URL}{\rm :} 
%http://www-math.mit.edu/\~{}lenny/}
%\renewcommand{\subjclassname}{\textup{2000} Mathematics Subject Classification}
%\subjclass{Primary 57R17; secondary 53D12, 57M27.}

\begin{abstract}
We establish tools to facilitate the computation and application
of the Chekanov-Eliashberg differential graded
algebra (DGA), a Legendrian-isotopy
invariant of Legendrian knots and links in standard contact three-space.
More specifically, we reformulate the DGA 
in terms of front projections, and introduce
the characteristic algebra, a new invariant
derived from the DGA.  We use our techniques to distinguish
between several previously indistinguishable Legendrian knots and links.
\end{abstract}

\maketitle

\vspace{24pt}

\section{Introduction}

A {\it Legendrian knot} in standard contact $\R^3$ is a knot 
which is 
everywhere tangent to the two-plane distribution induced by the contact
one-form $dz - y\,dx$.  Two Legendrian knots are {\it Legendrian isotopic}
if there is a smooth isotopy between them through Legendrian knots.

Broadly speaking, we wish to determine when two Legendrian knots
are Legendrian isotopic.
There are two ``classical'' invariants for knots
under Legendrian isotopy,
Thurston-Bennequin number $tb$ and rotation number $r$.
These form a complete set of invariants for some knots, including
the unknot \cite{bib:EF}, torus knots \cite{bib:EH}, 
and the figure eight knot \cite{bib:EH}.

However, there do exist non-isotopic
Legendrian knots of the same topological type with the same $tb$ and $r$.
The method for demonstrating this fact is a new, nonclassical invariant
independently introduced by Chekanov \cite{bib:Che} and 
Eliashberg \cite{bib:Eli}.
We will use Chekanov's combinatorial formulation of this
invariant, which we call the
Chekanov-Eliashberg differential graded algebra (DGA).
Chekanov introduced a concept of equivalence between DGAs which he
called stable tame isomorphism; then two Legendrian-isotopic
knots have equivalent DGAs.

The Chekanov-Eliashberg DGA was originally defined as an algebra 
over $\Z/2$ with grading over $\Z/(2r(K))$.  This has subsequently 
been lifted \cite{bib:ENS} to an algebra
over the ring $\Z[t,t^{-1}]$ with grading over $\Z$, 
by following the picture from
symplectic field theory \cite{bib:EGH}; this lifted algebra is what
we will actually call the Chekanov-Eliashberg DGA.

There are two standard methods to portray Legendrian
knots in standard contact $\R^3$, via projections to $\R^2$: the
{\it Lagrangian projection} to the $xy$ plane, and the
{\it front projection} to the $xz$ plane.  Chekanov and Eliashberg, 
motivated by the general framework of contact homology \cite{bib:Eli} and
symplectic field theory \cite{bib:EGH}, used the Lagrangian projection
in their setups.

When we attempt to apply the Chekanov-Eliashberg DGA to distinguish
between Legendrian knots, we encounter two problems.  The first is
that it is not easy to manipulate Lagrangian-projected knots.
Chekanov gives a criterion in \cite{bib:Che} for a knot diagram in
$\R^2$ to be the Lagrangian projection of a Legendrian knot, but it
remains highly nontrivial to determine by inspection when two Lagrangian 
projections represent Legendrian-isotopic knots.

\begin{figure}
\centerline{
\includegraphics[width=1.6in,angle=270]{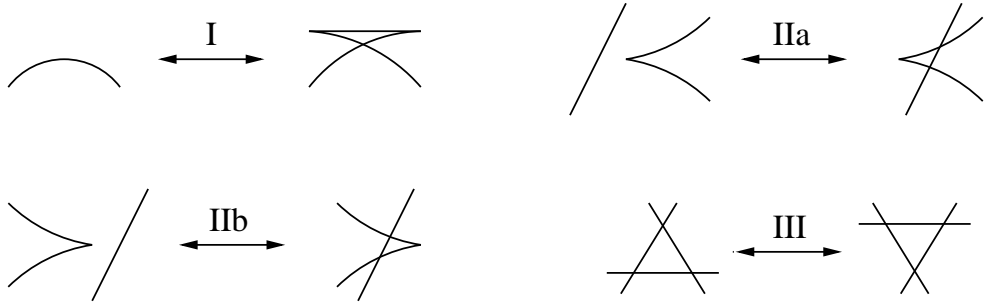}
}
\caption{The Legendrian Reidemeister moves which relate 
Legendrian-isotopic fronts.  Reflections of these moves in the horizontal
axis are also allowed.}
\label{fig:reidemeister}
\end{figure}

For questions of Legendrian
isotopy, the front projection is more convenient, because we know
precisely what diagrams represent front projections of Legendrian knots.
A {\it front} in $\R^2$ is simply a (continuous) embedding
of $S^1$ into the $xz$ plane, 
with a unique nonvertical tangent line
at each point in the image (i.e., so that $dz/dx$ exists at each
point), except of course at crossings.  Every front is the front
projection of a Legendrian knot, and every Legendrian knot projects
to such a front.  In addition, Legendrian-isotopic fronts are always
related by a series of Legendrian Reidemeister moves \cite{bib:Swi}:
see Figure~\ref{fig:reidemeister}.

The second problem is that it is difficult in general to tell
when two DGAs are equivalent.  To each DGA, Chekanov \cite{bib:Che} associated
an easy-to-compute Poincar\'e-type polynomial, which is invariant
under DGA equivalence, and used this to exhibit two $5_2$ knots which
have the same classical invariants but are not Legendrian isotopic.
On the other hand, the polynomial is only defined
for Legendrian knots possessing so-called augmentations;
in addition, there are often many nonisotopic knots with
identical classical invariants and polynomial.

This paper develops techniques designed to address
these problems.  In Section~\ref{sec:frontformulation},
we reformulate the Chekanov-Eliashberg DGA for front projections,
and discuss how it can often be easier to
compute than the Lagrangian-projection version.  
(The case of multi-component links, for which the results of this paper
hold with minor modifications, is explicitly addressed in
Section~\ref{ssec:dgalinks}.)
In 
Section~\ref{sec:charalg}, we introduce a new invariant, the characteristic
algebra, which is derived from the DGA and is 
relatively easy to compute.  The characteristic algebra is
quite effective in distinguishing between Legendrian isotopy classes;
it also encodes the information from both the Poincar\'e-Chekanov
polynomial and a similar higher-order invariant.  Section~\ref{sec:ex} 
applies the characteristic algebra to a number of knots and links which
were previously indistinguishable.

%*********************************************************************
\section{Chekanov-Eliashberg DGA in the front projection}
\label{sec:frontformulation}

This section is devoted to a reformulation of the Chekanov-Eliashberg
DGA from the Lagrangian projection to the more useful front projection.
In Section~\ref{ssec:morse}, we introduce resolution, the technique
used to translate from front projections to Lagrangian projections.
We then define the DGA for the front of a knot in
Section~\ref{ssec:front}, and discuss a particularly
nice and useful case in Section~\ref{ssec:simple}.  
In Section~\ref{ssec:invariance},
we review the main results concerning the DGA from \cite{bib:Che}
and \cite{bib:ENS}.  Section~\ref{ssec:dgalinks} discusses the
adjustments that need to be made for multi-component links.

%*********************************************************************
\subsection{Resolution of a front}
\label{ssec:morse}

Given a front, we can find a Lagrangian projection which represents
the same knot through the following construction, which is also considered
in \cite{bib:Fer} under the name ``morsification.''

\begin{figure}
\centering{
\includegraphics[height=5in,angle=270]{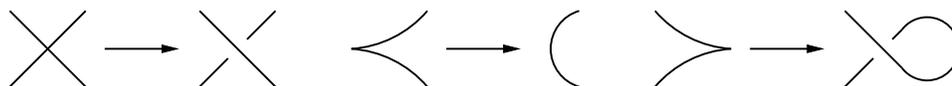}
}
\caption{
Resolving a front into the Lagrangian projection of a knot.
}
\label{fig:morse}
\end{figure}

\begin{definition}
The {\it resolution} of a front is the knot diagram obtained by
resolving each of the singularities in the front as shown in 
Figure~\ref{fig:morse}.
\label{not:morsification}
\end{definition}

\noindent 
The usefulness of this construction
is shown by the following result, which implies
that resolution is a map from front projections to Lagrangian
projections which preserves Legendrian isotopy.

\begin{proposition}
The resolution of the front projection of any Legendrian knot $K$
is the Lagrangian projection of a knot Legendrian
isotopic to $K$.
\label{prop:morse}
\end{proposition}

\noindent A similar result also holds for multi-component links.
Note that Proposition~\ref{prop:morse} 
is a bit stronger than the assertion
from \cite{bib:Fer} that the regular isotopy type of the resolution
is invariant under Legendrian isotopy of the front.

\begin{proof}
It suffices to distort the front $K$ smoothly to a front $K'$
so that the resolution of $K$ is the Lagrangian projection of 
the knot corresponding to $K'$.

\begin{figure}
\centering{
\includegraphics[height=3.6in,angle=270]{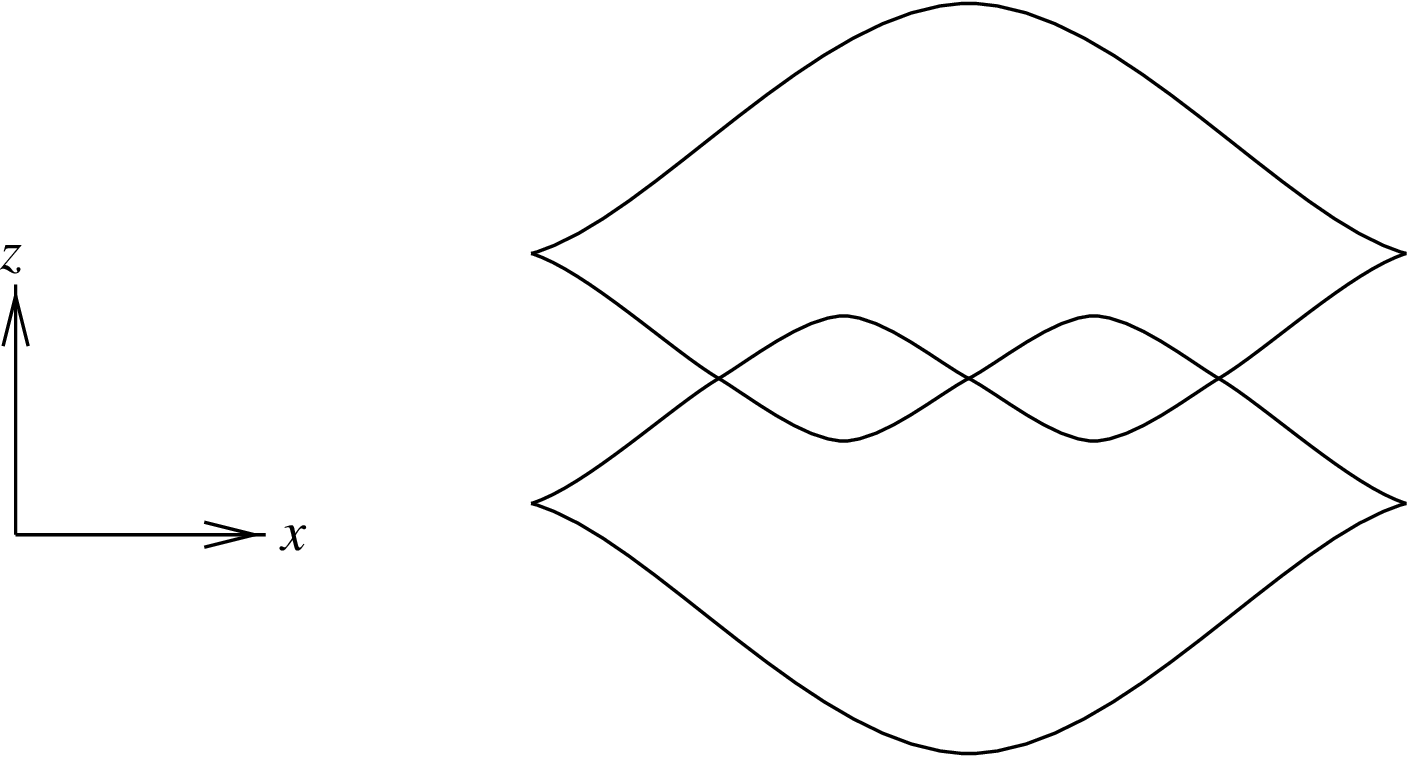}
\includegraphics[height=3.8in,angle=270]{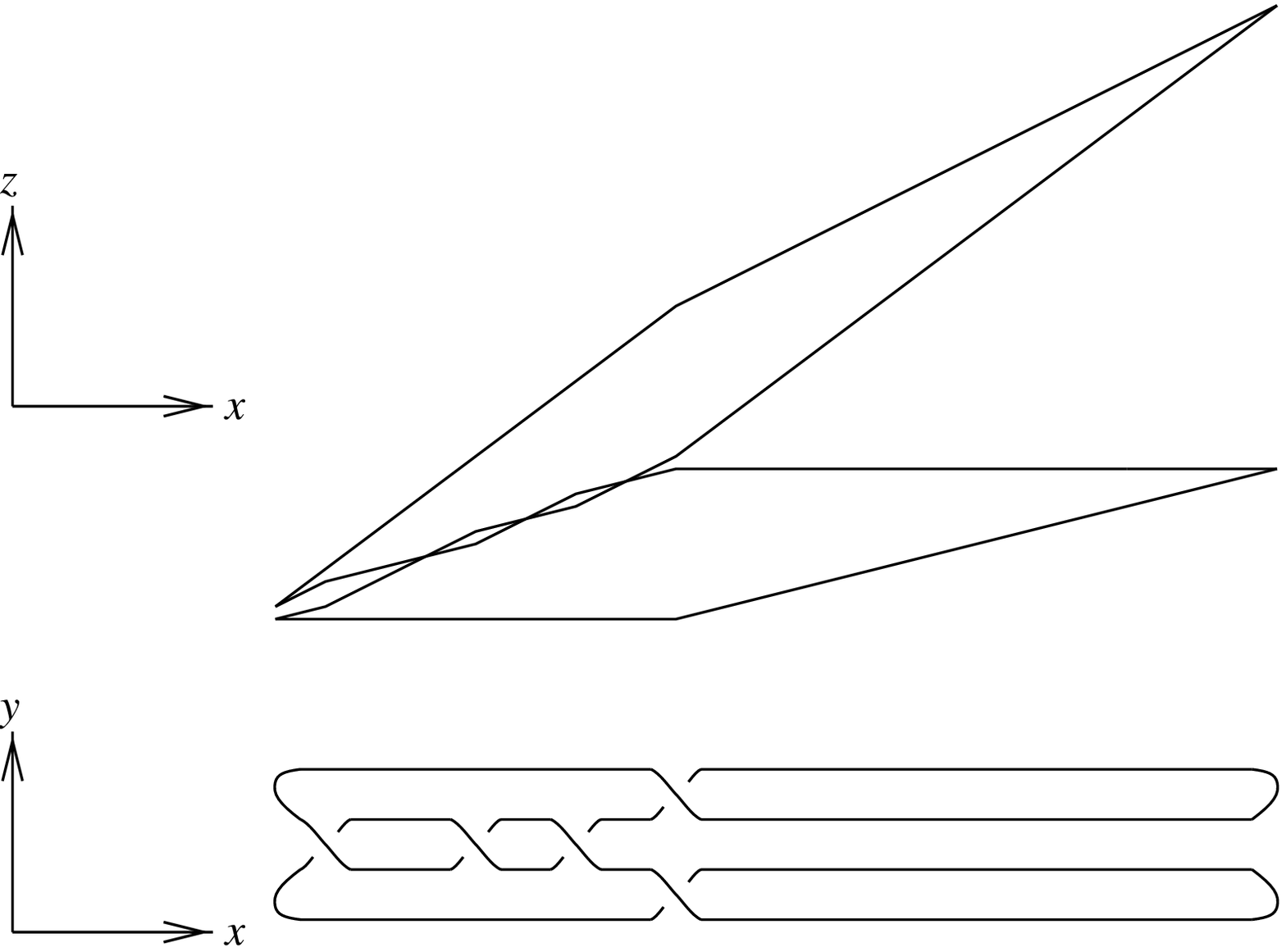}
}
\caption{
A front projection for the left-handed trefoil (top) is
distorted (middle) so that the corresponding Lagrangian projection (bottom),
given by $y=dz/dx$, with the same $x$ axis as the middle diagram,
is the resolution of the original front.  The exceptional segments
in the middle diagram appear as corners.
}
\label{fig:distort}
\end{figure}

We choose $K'$ to have the following 
properties; see Figure~\ref{fig:distort} for an illustration.  
Suppose that there are at most $k$ points in $K$ with any
given $x$ coordinate. 
Outside of arbitrarily small ``exceptional segments,'' $K'$ consists of 
straight line segments.  These line segments each have slope equal to
some integer between 0 and $k-1$ inclusive; outside of the exceptional
segments, for any given $x$ coordinate, the slopes of the line segments
at points with that $x$ coordinate are all distinct.  The purpose of
the exceptional segments is to allow the line segments to change slopes,
by interpolating between two slopes.
When two line segments exchange slopes via exceptional segments, the
line segment with higher $z$ coordinate has higher slope to the left
of the exceptional segment, and lower slope to the right.

It is always possible to construct such a distortion $K'$.  Build $K'$
starting from the left; a left cusp is simply two line segments of slope
$j$ and $j+1$ for some $j$, 
smoothly joined together by appending an exceptional segment to one
of the line segments.  Whenever two segments need to cross, force
them to do so by interchanging their slopes (again, with exceptional
segments added to preserve smoothness).  To create a right cusp between
two segments, interchange their slopes so that they cross, and then 
append an exceptional segment just before the crossing to preserve
smoothness.

We obtain the Lagrangian projection of the knot
corresponding to $K'$ by using the
relation $y=dz/dx$.  This projection consists of horizontal lines
(parallel to the $x$ axis), outside of a number of crossings arising from
the exceptional segments.  These crossings can be naturally identified
with the crossings and right cusps of $K$ or $K'$.  In particular,
right cusps in $K$ become the crossings associated to a simple loop.
It follows that the Lagrangian projection corresponding to $K'$
is indeed the resolution of $K$, as desired.
\end{proof}

%\begin{remark}
%Caveat lector: with the definitions we adopted in 
%Section~\ref{sec:intro2} for rotation numbers, the rotation number
%of the morsification of a front is {\it negative} the rotation number
%of the front.
%\end{remark}

%*********************************************************************
\subsection{The DGA for fronts of knots}
\label{ssec:front}

Suppose that we are given the front projection $Y$ of an oriented
Legendrian knot $K$.  To define the Chekanov-Eliashberg DGA
for $Y$, we simply examine the DGA for the resolution of $Y$
and ``translate'' this in terms of $Y$.  In the interests of readability,
we will concentrate on describing the DGA solely in terms of $Y$,
invoking the resolution only when the translation is not obvious.

\begin{figure}
\centering{
\includegraphics[width=1.8in,angle=270]{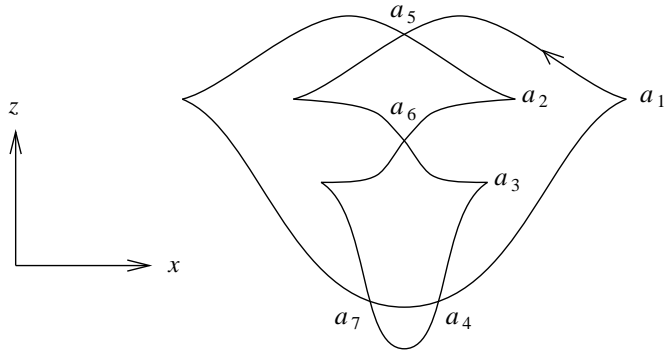}
}
\caption{The front projection of a figure eight knot,
with vertices labelled.}
\label{fig:figure8}
\end{figure}

The singularities of $Y$ fall into three
categories: crossings (nodes), left cusps, and right cusps.
Ignore the
left cusps, and call the crossings and right cusps {\it vertices},
with labels $a_1,\ldots,a_n$ (see Figure~\ref{fig:figure8});
then the vertices of $Y$ are in one-to-one correspondence with
the crossings of the resolution of $Y$.

As an algebra, the Chekanov-Eliashberg DGA of the front $Y$ is defined to be
the free, noncommutative, unital algebra $A=\alg{a_1,\ldots,a_n}$ 
over $\Z[t,t^{-1}]$ generated by $a_1,\ldots,a_n$.  We wish to
define a grading on $A$, and a differential $\d$ on $A$ which lowers the
grading by $1$.

We first address the grading of $A$.  For an oriented path 
$\gamma$ contained in the diagram $Y$, define $c(\gamma)$ 
\label{not:c} to be
the number of cusps traversed upwards, minus the number
of cusps traversed downwards, along $\gamma$.  
Note that this is the opposite convention from the one used to
calculate rotation number; if we consider $Y$ itself to be an
oriented closed curve, then $r(K) = -c(Y)/2$.

Let the degree of the indeterminate $t$ be $2r(K)$.
To grade $A$, it then suffices to define the degrees of the
generators $a_i$; we follow \cite{bib:ENS}.

\begin{definition}
Given a vertex $a_i$, define the {\it capping path}
$\gamma_i$, a path in $Y$ beginning and ending at $a_i$, as follows.
If $a_i$ is a crossing,
move initially along the segment of higher slope at $a_i$, in the direction
of the orientation of $Y$; then follow $Y$, not changing direction
at any crossing, until $a_i$ is reached again.  
If $a_i$ is a right cusp, then $\gamma_i$ is the
empty path, if the orientation of $Y$ traverses $a_i$ upwards,
or the entirety of $Y$ in the direction of its orientation,
if the orientation of $Y$ traverses $a_i$ downwards.
\label{def:cappingpath}
\end{definition}

\begin{definition}
If $a_i$ is a crossing, then $\deg a_i = c(\gamma_i)$.
If $a_i$ is a right cusp, then $\deg a_i$ is $1$ or $1 - 2r(K)$, depending
on whether the orientation of $Y$ traverses $a_i$ upwards or downwards,
respectively.
\label{def:degree}
\end{definition}

\noindent
We thus obtain a grading for $A$ over $\Z$.  
As an example, in the figure eight knot shown in Figure~\ref{fig:figure8},
$a_1,a_2,a_3,a_4,a_7$ have degree 1, while $a_5,a_6$ have degree 0.
%For an illustration of Definition~\ref{def:degree} for a knot
%of nonzero rotation number, see Remark~\ref{rmk:grading63}.

It will be useful
to introduce the sign function $\sgn v = (-1)^{\deg v}$
\label{not:sgn}
on pure-degree elements of $A$, including
vertices of $Y$; note that any right cusp has negative sign.
The Thurston-Bennequin number for $K$ can be written as the difference
between the numbers of positive-sign and negative-sign vertices in $Y$.
Since $\deg t = 2r(K)$,
the graded algebra $A$ incorporates both classical
Legendrian-isotopy invariants.

We next wish to define the differential $\d$ on $A$.  
As in \cite{bib:Che}, we define $\d a_i$ for a generator $a_i$
by considering a certain class of immersed disks in the diagram $Y$.

\begin{figure}
\centering{
\includegraphics[width=4in,angle=270]{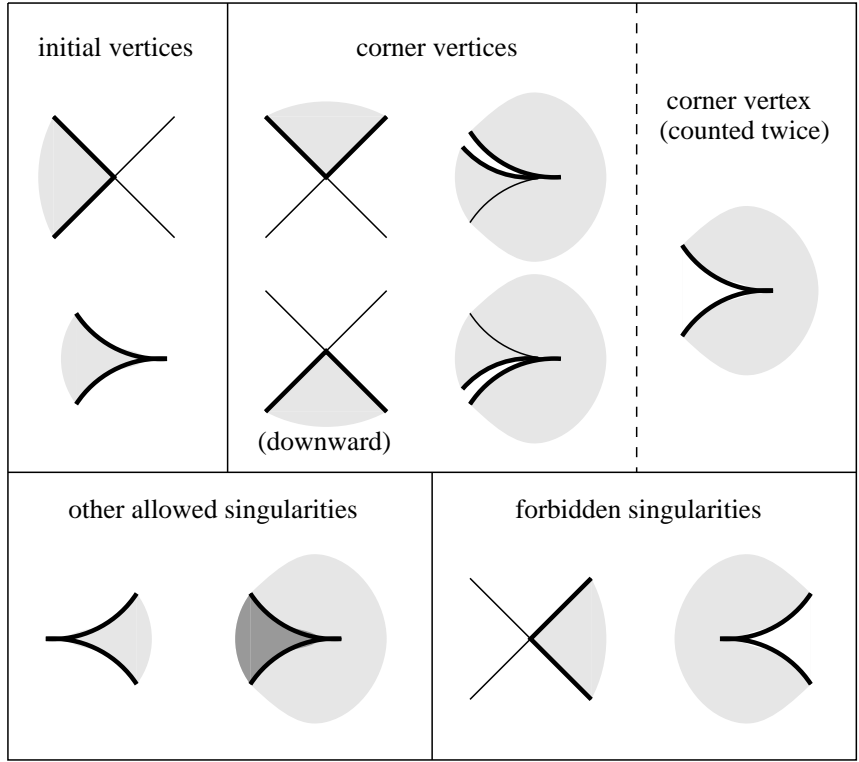}
}
\caption{
Possible singularities in an admissible map, and their classification.
The shaded area is the image
of the map restricted to a neighborhood of the singularity; the heavy
line indicates the image of the boundary of $D^2$.  In two of the diagrams,
the heavy line has been shifted off of itself for clarity.  The diagram with
heavy shading indicates that the image overlaps itself.  
%The last two
%diagrams are forbidden in an admissible map.
}
\label{fig:vertices}
\end{figure}

\begin{definition}
An {\it admissible map} on $Y$ is an immersion from the two-disk $D^2$ 
to $\R^2$ which maps
the boundary of $D^2$ into the knot projection $Y$, and which satisfies
the following properties: the map is smooth except possibly at
vertices and left cusps; the image of the map near any singularity 
looks locally like one of the diagrams in Figure~\ref{fig:vertices},
excepting the two forbidden ones; and,
in the notation of Figure~\ref{fig:vertices}, there is precisely one initial 
vertex.
\label{def:admissible}
\end{definition}

\noindent
The singularities of an admissible map thus consist
of one initial vertex, a number of corner vertices (possibly including
some right cusps counted twice), and some other singularities which we
will ignore.  One type of corner vertex, the ``downward'' corner vertex
as labelled in Figure~\ref{fig:vertices}, will be important shortly in
determining signs.

The possible singularities depicted in Figure~\ref{fig:vertices}
are all derived by considering the resolution of $Y$, but it is not
immediately obvious why the two forbidden singularities should be 
disallowed.  To justify this, call a point $p$ in the domain
of an admissible map, and its image under the map, 
{\it locally rightmost} if $p$ attains a local maximum for the $x$
coordinate of its image.
%  (More sloppily, a point in the image of
%the map is locally rightmost if it locally maximizes $x$ coordinate
%in the image.)  Observe that
Any locally rightmost point in the image
of an admissible map must be the unique initial vertex of the map:
this point must be a node or a right cusp, which cannot be a negative
corner vertex (cf.\ Figure~\ref{fig:vertices}).  In particular,
there must be a unique locally rightmost point in the image.  Of the two
forbidden singularities from Figure~\ref{fig:vertices}, the left one
is disallowed because the initial vertex is not rightmost, and the right
one because there would be two locally rightmost points.

To each diffeomorphism class of admissible maps on $Y$,
we will now associate a monomial in $\alg{a_1,\ldots,a_n}$.
Let $f$ be a representative of a diffeomorphism class, and suppose that
$f$ has corner vertices at $a_{j_1},\ldots,a_{j_{\ell}}$, counted twice
where necessary, in counterclockwise order around the boundary of $D^2$,
starting just after the initial vertex, and ending just before reaching 
the initial vertex again.  Then the monomial associated to $f$,
and by extension to the diffeomorphism class of $f$, is
\[
\alpha(f) = (\sgn f)\, t^{-n(f)} a_{j_1} \cdots a_{j_\ell},
\]
\label{not:alphaf}
where $(\sgn f)$ \label{not:sgnf} is the parity ($+1$ for even,
$-1$ for odd) of the number of
downward corner vertices of $f$ of even degree, and the 
winding number $n(f)$ is defined below.

The image $f(\d D^2)$, oriented counterclockwise, lifts to a collection
of oriented paths in the knot $K$.  If $a_i$ is the initial vertex of $f$,
then the lift of $f(\d D^2)$, along with the lifts of the
capping paths $\gamma_i$, $-\gamma_{j_1},\ldots,
-\gamma_{j_\ell}$, form a closed cycle in $K$.  We then set
$n(f)$ \label{not:nf} to be the winding number of this cycle 
around $K$, with respect to the orientation of $K$.

\begin{definition}
Given a generator $a_i$, we define
\[
\d a_i = 
\begin{cases}
\sum \alpha(f) & \text{if $a_i$ is a crossing} \\
1 + \sum \alpha(f) & \text{if $a_i$ is a right cusp oriented upwards} \\
t^{-1} + \sum \alpha(f) & \text{if $a_i$ is a right cusp oriented
downwards,}
\end{cases}
\]
where the sum is over all diffeomorphism classes of
admissible maps $f$ with initial vertex at $a_i$.
We extend the differential to the algebra $A$ by setting
$\d (\Z[t,t^{-1}])=0$ and imposing the signed
Leibniz rule $\d (vw) = (\d v) w + (\sgn v) v (\d w)$.
\label{def:diff}
\end{definition}

A few remarks are in order.  
The power of $t$ in the definition of the monomial $\alpha(f)$ has been
translated directly from the corresponding definition in \cite{bib:ENS}.  
It is easy to check that the signs also correspond 
to the signs in \cite{bib:ENS}, after we replace $a_i$ by $-a_i$ for each
$a_i$ which is ``right-pointing''; that is, near which the
knot is locally oriented from left to right for both strands.

Definition~\ref{def:diff} depends on a choice of orientation
of $K$.  For an unoriented knot, we may similarly define
the differential without the powers of $t$; the DGA is then an
algebra over $\Z$ graded over $\Z/(2r(K))$, still a lifting of
Chekanov's original DGA over $\Z/2$.

As a final remark, if $K$ is a stabilization, i.e., contains a zigzag
(see, e.g., \cite{bib:EH}), then it is easy to 
see that there is
an $a_i$ such that $\d a_i = 1$ or $\d a_i = t^{-1}$.  In this
case, $\d (a_j - a_i \d a_j) = 0$ or
$\d (a_j - t a_i \d a_j) = 0$ for all $j$, and the DGA
collapses modulo tame isomorphisms (see Section~\ref{ssec:invariance}).
This was first noted in \cite[\S 11.2]{bib:Che}.

\begin{figure}
\centering{
\includegraphics[height=5.5in,angle=270]{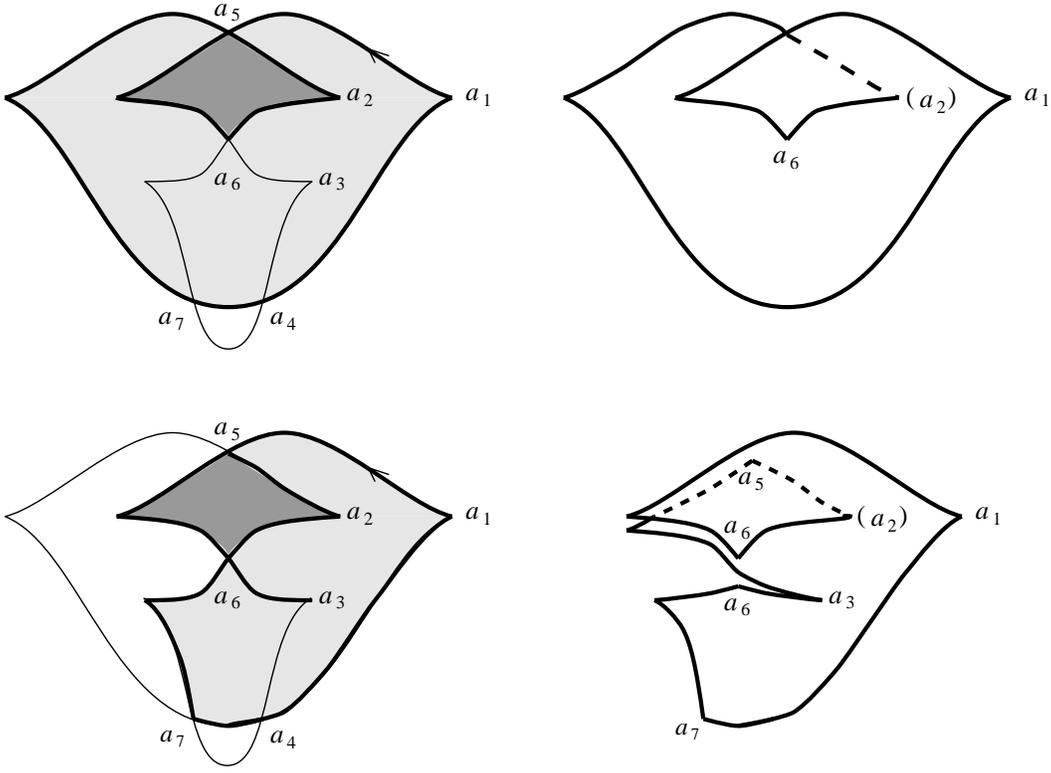}
}
\caption{
The admissible maps which correspond to the terms
$a_6$ (top) and $t^3 a_6a_5a_3a_6a_7$ (bottom) in $\d a_1$, for 
the front from
Figure~\ref{fig:figure8}.  The heavy lines indicate the image of the boundary
of $D^2$; the heavy shading indicates where the images overlap themselves.
For clarity, the images of the maps are redrawn to the right.
}
\label{fig:figure8gross}
\end{figure}

For the front in Figure~\ref{fig:figure8}, we may compute 
(somewhat laboriously) that
\begin{eqnarray*}
\d a_1 &=& 
1 + a_6 - t^2 a_6a_4a_6a_7 - t^2 (1-t a_6 a_5)a_3a_6a_7 
+ t a_6a_2 (1-ta_6-t^2 a_7a_4a_6)a_7 \\
\d a_2 &=& 1-ta_5a_6 \\
\d a_3 &=& t^{-1}-a_6-ta_6a_7a_4 \\
\d a_4 &=& \d a_5 = \d a_6 = \d a_7 = 0.
\end{eqnarray*}
See Figure~\ref{fig:figure8gross} for a depiction of two of the 
admissible maps counted in $\d a_1$.

To illustrate the calculation of the sign and power of $t$ associated to
an admissible map, consider the term $t^3 a_6 a_5 a_3 a_6 a_7$ in
$\d a_1$ above.  The sign of this term is $(-\sgn a_5)(-\sgn a_6) = +1$.
To calculate the power of $t$, we count, with orientation,
 the number of times the cycle corresponding to this map passes through 
$a_1$.  The boundary of the immersed disk passes through $a_1$, contributing
$1$; $\gamma_1$ trivially does not pass
through $a_1$, contributing $0$; and 
$-\gamma_3,-\gamma_6,-\gamma_7$ pass through $a_1$,
while $-\gamma_5$ does not, contributing a total of $-4$.  
It follows that the power of $t$ is $t^{-(1+0-4)} = t^3$.

%*********************************************************************
\subsection{Simple fronts}
\label{ssec:simple}

Since the behavior of an admissible map near a right cusp can be complicated, 
our formulation of the differential algebra may seem no easier to compute
than Chekanov's.  There is, however, one class of fronts for which
the differential is particularly easy to compute.

\begin{definition}
A front is {\it simple}
if it is smoothly isotopic to a front all of whose right cusps have the
same $x$ coordinate.
\label{def:simple}
\end{definition}

\noindent
Any front can be Legendrian-isotoped
to a simple front: 
``push'' all of the right cusps to the right until they share the
same $x$ coordinate.  (In the terminology of Figure~\ref{fig:reidemeister},
a series of IIb moves can turn any front into a simple front.)

For a simple front, the boundary of any admissible map
must begin at a node or right cusp (the initial vertex),
travel leftwards to a left cusp, and then travel rightwards again
to the initial vertex.  Outside of the initial vertex and the left cusp,
the boundary can only have very specific corner vertices: each corner
vertex must be a crossing, and, in a neighborhood of each of these nodes,
the image of the map must only occupy one of the four regions surrounding
the crossing.  In particular, the map is an embedding, not just an immersion.

\begin{figure}
\centering{
\includegraphics[width=1.8in,angle=270]{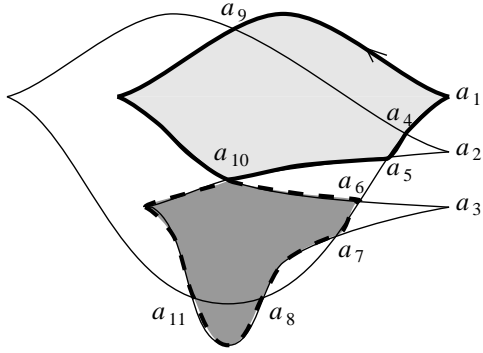}
}
\caption{
A simple-front version of the
front from Figure~\ref{fig:figure8}, with two admissible maps drawn.
The top shaded region corresponds to the term $ta_{10}a_5$ in $\d a_1$;
the bottom shaded region corresponds to the term $-ta_{10}a_7$ in $\d a_6$.
}
\label{fig:simple}
\end{figure}

\begin{example}
It is easy to calculate the differential for the
simple-front version of the figure eight knot given in
Figure~\ref{fig:simple}:
\[
\begin{array}{rclcrcl}
\d a_1 &=& 1+a_6+ta_{10}a_5 & \hspace{0.5in} & 
\d a_4 &=& t^{-1}+a_8a_7 - a_9a_6 - ta_9a_{10}a_5 \\
\d a_2 &=& 1-ta_9a_{10} && \d a_5 &=& a_7 + a_{11} + ta_{11}a_8a_7 \\
\d a_3 &=& t^{-1}-a_{10} - ta_{10}a_{11}a_8 && 
\d a_6 &=& - ta_{10}a_7 - ta_{10}a_{11} - t^2a_{10}a_{11}a_8a_7 \\
\multicolumn{7}{c}{
\d a_7 =\d a_8 =\d a_9 =\d a_{10} =\d a_{11} = 0.}
\end{array}
\]
For the signs, note that $a_1,a_2,a_3,a_4,$ and $a_8$ have degree 1,
$a_7$ and $a_{11}$ have degree $-1$, and the other vertices have degree 0;
for the powers of $t$, note that $\gamma_3$, $\gamma_4$, $\gamma_5$,
$\gamma_6$, $\gamma_7$, $\gamma_{10}$, and $\gamma_{11}$ pass through
$a_1$, while the other capping paths do not.
\end{example}

%*********************************************************************
\subsection{Properties of the DGA}
\label{ssec:invariance}

In this section, we summarize the properties of the Chekanov-Eliashberg
DGA.  These results were originally proven over $\Z/2$ in \cite{bib:Che},
and then extended over $\Z[t,t^{-1}]$ in \cite{bib:ENS}.
Proofs can be found in \cite{bib:ENS} in the Lagrangian-projection
setup, or in \cite{bib:thesis} in the front-projection setup.

\begin{proposition}[\cite{bib:Che},\cite{bib:ENS}]
For the DGA associated to a Legendrian knot, 
$\d$ lowers degree by $1$.
\label{prop:lowerdegree}
\end{proposition}

\begin{proposition}[\cite{bib:Che},\cite{bib:ENS}]
For the DGA associated to a Legendrian knot, $\d^2 = 0$.
\label{prop:dsquared}
\end{proposition}

To state the invariance result for the DGA, we need
to recall several definitions from \cite{bib:Che} or \cite{bib:ENS}.

An (algebra) automorphism of a graded free 
algebra $\alg{a_1,\ldots,a_n}$ is {\it elementary} \label{not:elementary}
if it preserves grading and
sends some $a_i$ to $a_i + v$, where $v$ does not involve $a_i$, and 
fixes the other generators $a_j, j\neq i$.  
A {\it tame automorphism} \label{not:tame} of $\alg{a_1,\ldots,a_n}$ is 
any composition of elementary automorphisms; a {\it tame isomorphism}
between two free algebras $\alg{a_1,\ldots,a_n}$ and $\alg{b_1,\ldots,b_n}$ 
is a grading-preserving 
composition of a tame automorphism and the map sending $a_i$ to $b_i$
for all $i$.  Two DGAs are then {\it tamely isomorphic} if there
is a tame isomorphism between them which maps the differential on one
to the differential on the other.

Let $E$ be a DGA with generators $e_1$ and $e_2$, such that
$\d e_1=\pm e_2$, $\d e_2=0$, both $e_1$ and $e_2$ have pure degree, and
$\deg e_1 = \deg e_2 + 1$.
Then an {\it algebraic stabilization} of a DGA 
$(A=\alg{a_1,\ldots,a_n},\d)$ is a graded coproduct
\[
(S(A),\d) = (A,\d) \amalg (E,\d) = (\alg{a_1,\ldots,a_n,e_1,e_2},\d),
\]
with differential and grading induced from $A$ and $E$.
Finally, two DGAs are {\it equivalent} if they are tamely isomorphic
after some (possibly different) number of (possibly different) 
algebraic stabilizations of each.

We can now state the main invariance result.

\begin{theorem}[\cite{bib:Che},\cite{bib:ENS}] \label{thm:invariance}
Fronts of Legendrian-isotopic knots have equivalent DGAs.
\end{theorem}

\begin{corollary}[\cite{bib:Che},\cite{bib:ENS}]
The graded homology of the DGA associated to a Legendrian knot is
invariant under Legendrian isotopy.
\label{cor:homology}
\end{corollary}

%*********************************************************************
\subsection{The DGA for fronts of links}
\label{ssec:dgalinks}

In this section, we describe the modifications of the definition
of the Chekanov-Eliashberg DGA necessary for Legendrian links
in standard contact $\R^3$.  Here the DGA has an infinite family of
gradings, as opposed to one, and is defined over a ring more complicated
than $\Z[t,t^{-1}]$.  The DGA for links also
includes some information not found for knots.

Let $L$ be an oriented Legendrian link, with 
components $L_1,\ldots,L_k$; in this section,
for ease of notation, we will also
use $L,L_1,\ldots,L_k$ to denote the corresponding front projections.
Chekanov's original definition \cite{bib:Che} of the DGA for $L$ gives
an algebra over $\Z/2$ graded over $\Z/(2r(L))$, where 
$r(L) = \gcd(r(L_1),\ldots,r(L_k))$; we will extend this to
an algebra over $\Z[t_1,t_1^{-1},\ldots,t_k,t_k^{-1}]$ graded
over $\Z$, and our set of gradings will be more refined than in
\cite{bib:Che}.  We will also discuss an additional structure on the DGA
introduced by K.\ Michatchev \cite{bib:Mi}.

As in Section~\ref{ssec:front}, let $a_1,\ldots,a_n$ be the vertices
(crossings and right cusps) of $L$.  
We associate to $L$ the algebra
\[
A = \Z[t_1,t_1^{-1},\ldots,t_k,t_k^{-1}] \langle a_1,\ldots,a_n \rangle,
\]
with differential and grading to be defined below.

For each crossing $a_i$, let $N_u(a_i)$ and $N_l(a_i)$ denote neighborhoods 
of $a_i$ on the two strands intersecting at $a_i$, so that
the slope of $N_l(a_i)$ is greater than the slope of $N_u(a_i)$, i.e.,
$N_u(a_i)$ is {\it lower} than $N_l(a_i)$ in $y$ coordinate.  
If $a_i$ is a right cusp, define $N_u(a_i) = N_l(a_i)$
to be a neighborhood of $a_i$ in $L$.
For any vertex $a_i$, we may then define two numbers $u(a_i)$ and $l(a_i)$,
\label{not:ul}
the indices of the link components containing $N_u(a_i)$ and $N_l(a_i)$,
respectively.

For each $j=1,\ldots,k$, fix a base point $p_j$ on $L_j$, away from
the singularities of $L$, so that
$L_j$ is oriented from left to right in a neighborhood of $p_j$.
To a crossing $a_i$, we associate two
capping paths $\gamma^u_i$ and $\gamma^l_i$: \label{not:gammaul}
$\gamma^u_i$ is the path beginning at
$p_{u(a_i)}$ and following $L_{u(a_i)}$ in the direction of its
orientation until $a_i$ is reached through $N_u(a_i)$; 
$\gamma^l_i$ is the analogous path
in $L_{l(a_i)}$ beginning at $p_{l(a_i)}$ and ending at $a_i$ through
$N_l(a_i)$.  (If $u(a_i) = l(a_i)$, then one of $\gamma^u_i$ and 
$\gamma^l_i$ will contain the other.)
Note that, by this definition, when $a_i$ is a right cusp, 
$\gamma^u_i$ and $\gamma^l_i$ are both the path beginning at
$p_{u(a_i)} = p_{l(a_i)}$ and ending at $a_i$.

\begin{definition}
For $(\rho_1,\ldots,\rho_{k-1}) \in \Z^{k-1}$, we may define
a $\Z$ grading on $A$ by
\[
\deg a_i = \begin{cases}
1 & \text{if $a_i$ is a right cusp} \\
c(\gamma^u_i) - c(\gamma^l_i) + 2\rho_{u(a_i)} - 2\rho_{l(a_i)}
& \text{if $a_i$ is a crossing,}
\end{cases}
\]
where we set $\rho_k = 0$.  We will only consider gradings on $A$
obtained in this way.
\label{def:linkdegree}
\end{definition}

\noindent
The set of gradings on $A$ is then indexed by $\Z^{k-1}$.
(In particular, a knot has precisely one grading, the one given
in Section~\ref{ssec:front}.)
Our motivation for including precisely this set of gradings is given
by the following easily proven observation.

\begin{lemma}
The collection of possible gradings on $A$
is independent of the choices of the points $p_j$.
\end{lemma}

If $a_i$ is contained entirely in component $L_j$, then the degree of $a_i$ 
may differ from how we defined it in Definition~\ref{def:degree}
with $L_j$ a knot by itself.
It is easy to calculate that the difference between the two degrees 
will always be either $0$ or $2r(L_j)$.

We may define the sign function on vertices, as usual, by
$\sgn a_i = (-1)^{\deg a_i}$.
This is well-defined and independent of the choice of grading:
$\sgn a_i = -1$ if $a_i$ is a right cusp; $\sgn a_i = 1$ if
$a_i$ is a crossing with both strands pointed in the same direction
(either both to the left or both to the right); and $\sgn a_i = -1$
if $a_i$ is a crossing with strands pointed in opposite directions.
Note that $tb(L) = \sum_{i=1}^n \sgn a_i$.

We may still define the differential of a generator $a_i$ as in
Definition~\ref{def:diff}, but we must now redefine 
$\alpha(f)$ for an admissible map $f$.  Suppose that $f$ has initial
vertex $a_i$ and corner vertices $a_{i_1},\ldots,a_{i_m}$.
Then the lift of $f(\d D^2)$ to $L$,
together with the lifts of $\gamma^u_i,$
$-\gamma^u_{i_1},\ldots,-\gamma^u_{i_m}$, $-\gamma^l_i,$ 
$\gamma^l_{i_1},\ldots,\gamma^l_{i_m},$
form a closed cycle in $L$.
Let the winding number of this cycle around component $L_j$ be
$n_j(f)$.  Also, define $\sgn f$,
as before, to be the parity of the number of downward corner vertices
of $f$ with positive sign.  

We now set
\[
\alpha(f) = (\sgn f) \, t_1^{-n_1(f)} \cdots t_k^{-n_k(f)}
a_{i_1} \cdots a_{i_m}.
\]
The differential $\d$ can then be defined on $A$ essentially as in 
Definition~\ref{def:diff}, except that we now have
$\d(\Z[t_1,t_1^{-1},\ldots,t_k,t_k^{-1}]) = 0$, and
\[
\d a_i = 
\begin{cases}
\sum \alpha(f) & \text{if $a_i$ is a crossing} \\
1 + \sum \alpha(f) & 
\text{if $a_i$ is a right cusp.}
\end{cases}
\]
Note that the signed Leibniz rule does not depend on the choice
of base points $p_j$, since
the signs $(\sgn a_i)$ are independent of this choice.  
Also, because of a different choice of
capping paths, we always add $1$ to a right cusp; cf.\ 
Definition~\ref{def:diff}.

In practice, there is a simple way to calculate $n_j(f)$: 
it is the signed number
of times $f(\d D^2)$ crosses $p_j$.  Indeed, the winding number
of the appropriate cycle around $L_j$ is the signed number of times that
it crosses a point on $L_j$ just to the left of $p_j$.   No
capping path $\gamma_i^u$ or $\gamma_i^l$, however, crosses this point.
Hence $n_j(f)$ counts the number of times $f(\d D^2)$ crosses a point just
to the left of $p_j$; we could just as well consider $p_j$ instead of
this point.

We next examine the effect
of changing the base points $p_j$ on the differential $\d$.  Consider
another set of base points $\tilde{p}_j$, giving rise to capping paths
$\tilde{\gamma}^u_i,\tilde{\gamma}^l_i$, and let $\xi_j$ be the
oriented path in $L_j$ from $p_j'$ to $p_j$.  Then
\[
\tilde{\gamma}^u_i - \gamma^u_i = 
\begin{cases}
\xi_{u(a_i)}, & N_u(a_i) \subset \xi_{u(a_i)} \\
\xi_{u(a_i)} - L_{u(a_i)}, & N_u(a_i) \not\subset \xi_{u(a_i)},
\end{cases}
\]
and similarly for $\tilde{\gamma}^l_i - \gamma^l_i$.
We conclude the following result.

\begin{lemma} \label{lem:baseptchange}
The differential on $A$, calculated with base points
$\tilde{p}_j$, is related to the differential calculated with $p_j$,
by intertwining with the following automorphism on $A$:
\[
a_i \mapsto \begin{cases}
a_i, & \text{$N_u(a_i) \subset \xi_{u(a_i)}$ and 
$N_l(a_i) \subset \xi_{l(a_i)}$} \\
t_{l(a_i)}^{-1} a_i, & \text{$N_u(a_i) \subset \xi_{u(a_i)}$ and 
$N_l(a_i) \not\subset \xi_{l(a_i)}$} \\
t_{u(a_i)} a_i, & \text{$N_u(a_i) \not\subset \xi_{u(a_i)}$ and 
$N_l(a_i) \subset \xi_{l(a_i)}$} \\
t_{u(a_i)} t_{l(a_i)}^{-1} a_i, & \text{$N_u(a_i) \not\subset 
\xi_{u(a_i)}$ and $N_l(a_i) \not\subset \xi_{l(a_i)}$.}
\end{cases}
\]
\end{lemma}

\begin{figure}
\centering{
\includegraphics[width=1.5in,angle=270]{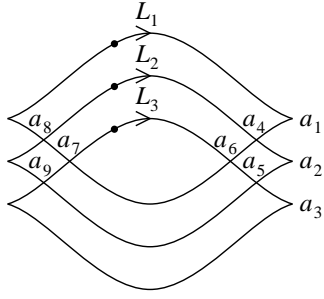}
}
\caption{
An oriented link $L$ with components $L_1$, $L_2$, and
$L_3$, with corresponding base points $p_1$, $p_2$, and
$p_3$ marked but not labelled.
}
\label{fig:triple}
\end{figure}

%\begin{example}
%\label{ex:triple}
For example, consider the link $L$ in Figure~\ref{fig:triple}, with base
points as shown.  To give a grading to the DGA on $L$, choose
$(\rho_1,\rho_2) \in \Z^2$.
The degrees of generators are as follows:
\begin{alignat*}{3}
\deg a_1 &= 1 & \qquad \deg a_4 &= 1 + 2\rho_2 - 2\rho_1 &
\qquad \deg a_7 &= -1 + 2\rho_1 \\
\deg a_2 &= 1 & \deg a_5 &= 1 - 2\rho_2 &
\deg a_8 &= -1 + 2\rho_1 - 2\rho_2 \\
\deg a_3 &= 1 & \deg a_6 &= 1 - 2\rho_1 &
\deg a_9 &= -1 + 2\rho_2.
\end{alignat*}
The differential $\d$ is then given by
\begin{alignat*}{3}
\d a_1 &= 1 + t_1 + t_1t_2^{-1}a_8a_4 + t_1t_3^{-1}a_7a_6 & \qquad
\d a_4 &= t_2t_3^{-1}a_9 a_6 & \qquad \d a_7 &= a_8a_9 \\
\d a_2 &= 1 + t_2 + t_2t_3^{-1}a_9a_5 + a_4a_8 &
\d a_5 &= a_6a_8 & \d a_8 &= 0 \\
\d a_3 &= 1 + t_3 + a_5a_9 + a_6a_7 &
\d a_6 &= 0 & \d a_9 &= 0.
\end{alignat*}
%\end{example}

We can now state several properties of the link DGA, the analogues
of the results for knots in Section~\ref{ssec:invariance}.

\begin{proposition}
If $(A,\d)$ is a DGA associated to the link $L$, then $\d^2 = 0$, and
$\d$ lowers degree by $1$ for any of the gradings of $A$.
\label{prop:dsquaredlink}
\end{proposition}

The main invariance result requires a slight tweaking of the definitions.
Define elementary and tame automorphisms as in 
Section~\ref{ssec:invariance}; now, however, let a tame isomorphism
between algebras generated by $a_1,\ldots,a_n$ and $b_1,\ldots,b_n$
be a grading-preserving composition of a tame automorphism and
a map sending $a_i$ to 
$\left( \prod_{m=1}^k t_m^{\nu_{k,m}} \right) b_i$, for
any set of integers $\{\nu_{k,m}\}$.
(This definition is necessitated by Lemma~\ref{lem:baseptchange}.)
Define algebraic stabilization and equivalence as before.

\begin{proposition}
If $L$ and $L'$ are Legendrian-isotopic oriented links,
then for any grading of the DGA for $L$, there is a grading
of the DGA for $L'$ so that the two DGAs are equivalent.
\label{prop:invariancelink}
\end{proposition}
\noindent
The proofs of Propositions~\ref{prop:dsquaredlink} and
\ref{prop:invariancelink} will be omitted here, as they are 
simply variants on the proofs of
Propositions \ref{prop:lowerdegree} and \ref{prop:dsquared} and
Theorem~\ref{thm:invariance}; see \cite{bib:Che} or \cite{bib:thesis}.

Our set of gradings for $A$ is more restrictive than the set
of ``admissible gradings'' postulated in \cite{bib:Che}.  To see this,
we first translate our criteria for gradings to the Lagrangian-projection
picture, and then compare with Chekanov's original criteria.

Consider a Legendrian link $L$ with components $L_1,\ldots,L_k$.  
By perturbing $L$
slightly, we may assume that the crossings of $\pixy(L)$ are
orthogonal, where $\pixy$ is the projection map $(x,y,z) \mapsto (x,y)$; 
as usual, label these crossings $a_1,\ldots,a_n$.  Choose neighborhoods
$N_u(a_i)$ and $N_l(a_i)$ in $L$ of the two points mapping to $a_i$
under $\pixy$, so that $N_u(a_i)$ lies above $N_l(a_i)$ in 
$z$ coordinate, and let $u(a_i)$ and $l(a_i)$ be the indices of
the link components on which these neighborhoods lie.

For each $j$, choose a point $p_j$ on $L_j$, and let
$\theta_j$ be an angle, measured counterclockwise, 
from the positive $x$ axis to the oriented tangent
to $L_j$ at $p_j$; note that $\theta_j$ is only well-defined up to
multiples of $2\pi$.
Let $r_u(a_i)$ be
the counterclockwise rotation number (the number of revolutions made)
for the path in $\pi(L_{u(a_i)})$ beginning at $p_{u(a_i)}$ and
following the orientation of $L_{u(a_i)}$ until $a_i$ is reached
via $N_u(a_i)$; similarly define $r_l(a_i)$.  Then the gradings for
the DGA of $L$ are given by choosing $(\rho_1,\ldots,\rho_{k-1})
\in \Z^{k-1}$ and setting
\[
\deg a_i = 2(r_l(a_i)-r_u(a_i))+(\theta_{l(a_i)}-\theta_{u(a_i)})/\pi
+ 2\rho_{u(a_i)} - 2\rho_{l(a_i)} - 1/2.
\]
By comparison, the allowed degrees in \cite{bib:Che} are given by
\[
\deg a_i = 2(r_l(a_i)-r_u(a_i))+(\theta_{l(a_i)}-\theta_{u(a_i)})/\pi
+ \rho_{u(a_i)} - \rho_{l(a_i)} - 1/2.
\]
The difference arises from the fact that Chekanov never uses the orientations
of the link components; this forces $\theta_{l(a_i)}$ and 
$\theta_{r(a_i)}$ to be well-defined only up to integer multiples of $\pi$,
rather than $2\pi$.

We now discuss an additional structure on the DGA for a link $L$, 
inspired by \cite{bib:Mi}.  More precisely, we will describe
a variant of the relative homotopy splitting from \cite{bib:Mi};
our variant will split something which is essentially a 
submodule of the DGA into $k^2$ pieces
which are invariant under Legendrian isotopy.

\begin{definition}
For $j_1 \neq j_2$ between $1$ and $k$, inclusive, define
$\Gamma_{j_1j_2}$ to be the module over 
$\Z[t_1,t_1^{-1},\ldots,t_k,t_k^{-1}]$ generated by words of the form
$a_{i_1} \cdots a_{i_m}$, with $u(a_{i_1}) = j_1,$
$l(a_{i_m}) = j_2$, and $u(a_{i_{p+1}}) = l(a_{i_p})$ for 
$1 \leq p \leq m-1$.
If $j_1 = j_2 = j$, then let $\Gamma_{j_1j_2}$ be the module generated by
such words, along with an indeterminate $e_j$.
Finally, let $\Gamma = \oplus \Gamma_{j_1j_2}$.
\end{definition}

\noindent
The indeterminates $e_j$ will replace the $1$ terms in the definition of
$\d$; see below.  Note that $a_i \in \Gamma_{u(a_i)l(a_i)}$.
Although $\Gamma$ itself is not an algebra, we have the usual
multiplication map
$\Gamma_{j_1j_2} \times \Gamma_{j_2j_3} \rightarrow \Gamma_{j_1j_3}$,
given on generators by concatenation, once we stipulate that the $e_j$'s
act as the identity.

Our introduction of $\Gamma$ is motivated by the fact that
$\d a_i$ is essentially in $\Gamma_{u(a_i)l(a_i)}$ for all $i$.  Define 
$\d' a_i$ as follows: if $u(a_i) \neq l(a_i)$, then 
$\d' a_i = \d a_i$; if $u(a_i) = l(a_i)$, then $\d' a_i$ is
$\d a_i$, except that we replace any $1$ or $2$ term in $\d a_i$
by $e_{u(a_i)}$ or $2e_{u(a_i)}$.  (It is easy to see that these
are the only possible terms in $\d a_i$ which involve only the $t_j$'s and no
$a_m$'s.)

\begin{lemma}
$\d' a_i \in \Gamma_{u(a_i)l(a_i)}$ for all $i$.
\label{lem:linkmodule}
\end{lemma}

\begin{proof}
For a term in $\d a_i$ of the form $a_{i_1}\ldots a_{i_k}$,
where we exclude powers of $t_j$'s,
we wish to prove that $u(a_{i_1}) = u(a_i)$, $l(a_{i_k}) = l(a_i)$,
and $u(a_{i_{p+1}}) = l(a_{i_p})$ for all $p$.  Consider the boundary
of the map which gives the term $a_{i_1}\ldots a_{i_k}$.  By definition,
the portion of this boundary connecting $a_{i_p}$ to $a_{i_{p+1}}$
belongs to link component
$l(a_{i_p})$ on one hand, and $u(a_{i_{p+1}})$ on the other.
We similarly find that $u(a_{i_1}) = u(a_i)$ and $l(a_{i_k}) = l(a_i)$.
\end{proof}

\begin{definition}
The {\it differential link module} of $L$ is $(\Gamma,\d')$, where we have
defined $\d' a_i$ above, and we
extend $\d'$ to $\Gamma$ by applying the signed Leibniz rule
and setting $\d' e_j = 0$ for all $j$.  A grading for $\Gamma$
is one inherited from the DGA of $L$, with $\deg e_j = 0$ for all $j$.
\label{def:linkmodule}
\end{definition}

We may define (grading-preserving) 
elementary and tame automorphisms and tame isomorphisms
for differential link modules as for DGAs, with the additional 
stipulation that all maps must preserve the link module structure 
by preserving $\Gamma_{j_1j_2}$ for all $j_1,j_2$.  Similarly, 
we may define an algebraic stabilization of a differential link module, with
the additional stipulation that the two added generators both
belong to the same $\Gamma_{j_1j_2}$.  As usual, we then define two
differential link modules to be equivalent if they are tamely
isomorphic after some number of algebraic stabilizations.
We omit the proof of the following result, which again is
simply a variant on the proof of the corresponding result for knots.

\begin{proposition}
If $L$ and $L'$ are Legendrian-isotopic oriented links,
then for any grading of the differential link module for $L$, there
is a grading of the differential link module for $L'$ so that
the two are equivalent.
\label{prop:linkmodule}
\end{proposition}

In this paper, we will not use the full strength of the
differential link module.
%, although applications of the characteristic
%module, which is derived from the differential link module (see
%Remark~\ref{rmk:charmodule}), are given in \cite{bib:NT} for links
%on the solid torus.  
We will, however, apply first-order
Poincar\'e-Chekanov polynomials
derived from the differential link module; we now describe these
polynomials, first mentioned in \cite{bib:Mi}.
For the definition of augmentations for knots, and background on 
Poincar\'e-Chekanov polynomials, please refer to 
\cite{bib:Che} or Section~\ref{ssec:polynomial}.

Assume that $r(L_1) = \cdots = r(L_k) = 0$, and let $\Gamma$ be
the differential link module for $L$, with some fixed grading.
We consider the DGAs for $L$ and $L_1,\ldots,L_k$ over $\Z/2$; that is,
set $t_j = 1$ for all $j$, and reduce modulo 2.

\begin{definition}
Suppose that, when considered alone as a knot, the DGA for each of 
$L_1,\ldots,L_k$
has an augmentation $\eps_1,\ldots,\eps_k$.  Extend these augmentations
to all vertices $a_i$ of $L$ by setting
\[
\eps(a_i) = \begin{cases}
\eps_{u(a_i)}(a_i) & \text{if $u(a_i) = l(a_i)$} \\
0 & \text{otherwise.}
\end{cases}
\]
We define an {\it augmentation} of $L$ to be any function $\eps$ 
obtained in this way.
\end{definition}

An augmentation $\eps$, as usual, gives rise to a first-order
Poincar\'e-Chekanov polynomial $P^{\eps,1}(\lambda)$; 
we may say, a bit imprecisely, that this polynomial splits into
$k^2$ polynomials $P^{\eps,1}_{j_1j_2}(\lambda)$, corresponding
to the pieces in $\Gamma_{j_1j_2}$.

The polynomials
$P^{\eps,1}_{jj}(\lambda)$ are precisely the polynomials
$P^{\eps_j,1}(\lambda)$ for each individual link component $L_j$.
For practical purposes, we can define $P^{\eps,1}_{j_1j_2}(\lambda)$ for
$j_1 \neq j_2$ as follows.  For $a_i \in \Gamma_{j_1j_2}$, define
$\d_\eps^{(1)} a_i$ to be the image of $\d a_i$ under the following
operation: discard all terms in $\d a_i$ containing more than one
$a_m$ with $u(a_m) \neq l(a_m)$, and 
replace each $a_m$ in $\d a_i$ by $\eps(a_m)$ whenever $u(a_m) = l(a_m)$.  
If we write $V_{j_1j_2}$ as the vector space over $\Z/2$
generated by $\{a_i \in \Gamma_{j_1j_2}\}$, then $\d_\eps^{(1)}$
preserves $V_{j_1j_2}$ and $\left(\d_\eps^{(1)}\right)^2 = 0$.
We may then set $P^{\eps,1}_{j_1j_2}(\lambda)$ to be the
Poincar\'e polynomial of $\d_\eps^{(1)}$ on $V_{j_1j_2}$, i.e.,
the polynomial in $\lambda$ whose $\lambda^i$ coefficient is
the dimension of the $i$-th graded piece of
$(\ker \d_\eps^{(1)}) / (\im \d_\eps^{(1)})$.

We may also define higher-order Poincar\'e-Chekanov polynomials
$P^{\eps,n}_{j_1j_2}(\lambda)$ by examining the action of $\d'$
on $\Gamma_{j_1j_2}$, but we will not need these here.

The following result,
which follows directly from Proposition~\ref{prop:linkmodule} and Chekanov's
corresponding result from \cite{bib:Che},
will be used in Section~\ref{sec:ex}.  

\begin{theorem} \label{thm:polylink}
Suppose that $L$ and $L'$ are Legendrian-isotopic oriented links.
Then, for any given grading and augmentation of the DGA for $L$, 
there is a grading and augmentation of the DGA for $L'$ so that
the first-order Poincar\'e-Chekanov polynomials $P^{\eps,1}_{j_1j_2}$
for $L$ and $L'$ are equal {\it for all} $j_1,j_2$.
\end{theorem}

For unoriented links,
we simply expand the set of allowed gradings $(\rho_1,\ldots,\rho_{k-1})$
to allow half-integers, as in \cite{bib:Che}.
Indeed, a grading of half-integers $(\rho_1,\ldots,\rho_{k-1})$
corresponds to changing the original orientation of $L$ by either
reversing the orientation of $\{L_j\,:\,2\rho_j~\mbox{odd}\}$,
or reversing the orientations of $L_k$ and
$\{L_j\,:\,2\rho_j~\mbox{even}\}$.  
We deduce this by examining how the
capping paths and degrees change when we change the orientation
(and hence base point) of one link component $L_j$.

%\begin{remark}
%While $P^{\eps,1}_{jj}(-1) = tb(L_j)$ as usual, we also have
%$P^{\eps,1}_{j_1j_2}(-1) = \operatorname{lk}(L_{j_1},L_{j_2})$, the
%linking number of $L_{j_1}$ and $L_{j_2}$, for $j_1 \neq j_2$.
%Also, we have $\sum_{j_1,j_2} P^{\eps,1}_{j_1j_2}(-1) = tb(L)$.
%We conclude that the first-order Poincar\'e-Chekanov polynomials 
%incorporate the classical invariants for oriented links (see
%Section~\ref{sec:intro2}).
%\end{remark}

%The following result demonstrates the effect of changing the
%orientation of one of the components on the polynomials $P^{\eps,1}_{j_1j_2}$.
%We will use this lemma in Chapter~\ref{chap:ex}.
%
%\begin{lemma}
%Suppose that we change the orientation (and hence base point) on
%one link component $L_j$, while keeping the orientations and base points
%fixed on all other link components.  Then, for all $j_1,j_2 \neq j$, 
%the polynomials
%$P^{\eps,1}_{jj_2}$ and $P^{\eps,1}_{j_1j}$ change by multiplication
%by $t_j^{2m+1}$ and $t_j^{-2m-1}$ for some fixed $m$.
%\end{lemma}
%
%\begin{proof}
%If $\gamma$ is either path in $L_j$ connecting the old base point to
%the new base point, then write $c(\gamma) = 2m+1$ with respect to
%the old orientation.  It is easy to see that the degrees of vertices
%yuck.

%\begin{example}
For the link from Figure~\ref{fig:triple}, an augmentation is any map with 
$\eps(a_i) = 0$ for $i \geq 4$.
Then $\d_\eps^{(1)}$ is identically zero, and the first-order
Poincar\'e-Chekanov polynomials simply measure the degrees of the 
$a_i$.  More precisely, for a choice of grading $(\rho_1,\rho_2) \in \Z^2$,
we have
\begin{alignat*}{3}
P^{\eps,1}_{11}(\lambda) &= \lambda & \qquad 
P^{\eps,1}_{21}(\lambda) &= \lambda^{1 + 2\rho_2 - 2\rho_1} & \qquad
P^{\eps,1}_{31}(\lambda) &= \lambda^{1 - 2\rho_1} \\
P^{\eps,1}_{12}(\lambda) &= \lambda^{-1 + 2\rho_1 - 2\rho_2} & \qquad 
P^{\eps,1}_{22}(\lambda) &= \lambda & \qquad
P^{\eps,1}_{32}(\lambda) &= \lambda^{1 - 2\rho_2} \\
P^{\eps,1}_{13}(\lambda) &= \lambda^{-1 + 2\rho_1} & \qquad 
P^{\eps,1}_{23}(\lambda) &= \lambda^{-1 + 2\rho_2} & \qquad
P^{\eps,1}_{33}(\lambda) &= \lambda.
\end{alignat*}
\label{ex:polytriple}
%\end{example}

%\begin{definition}
%An {\it augmentation} of $\Gamma$ is a map 
%$\eps\,:\,\Gamma \rightarrow \Z/2$ satisfying the following properties:
%\begin{enumerate}
%\item $\eps(t_j) = \eps(t_j^{-1}) = \eps(e_j) = 1$ for all $j$;
%\item $\eps(a_i) = 0$ if $u(a_i) \neq l(a_i)$ or $\deg a_i \neq 0$;
%\item $\eps$ commutes with multiplication on $\Gamma$ where defined;
%\item $\eps \circ \d = 0$.
%\end{enumerate}
%\end{definition}
%
%\noindent In other words, an augmentation of the differential link module
%is derived from an augmentation of the DGA which preserves the link
%module structure.  It is easy to see that an augmentation of $\Gamma$ is
%given precisely by augmentations on each of the DGAs of the link components.

\vspace{24pt}

%*********************************************************************

\section{The characteristic algebra}
\label{sec:charalg}

%\section{Introduction}

We would like to use the Chekanov-Eliashberg DGA to distinguish
between Legendrian isotopy classes of knots.  Unfortunately, it is often
hard to tell when two DGAs are equivalent.  In particular, the
homology of a DGA is generally infinite-dimensional and difficult to
grasp; this prevents us from applying Corollary~\ref{cor:homology}
directly.

Until now, the only known ``computable'' Legendrian invariants---that is,
nonclassical invariants which can be used in practice to distinguish between
Legendrian isotopy classes of knots---were the first-order Poincar\'e-Chekanov
polynomial and its higher-order analogues.  
However, the Poincar\'e-Chekanov polynomial is not defined 
for all Legendrian knots,
nor is it necessarily uniquely defined; in addition, as we shall see,
there are many nonisotopic knots with the same
polynomial.  The higher-order polynomials, 
on the other hand, are difficult to compute, and have not yet been 
successfully used to distinguish Legendrian knots.

In Section~\ref{ssec:charalgdef}, we introduce the characteristic algebra,
a Legendrian invariant derived from the DGA, which is nontrivial
for most, if not all, Legendrian knots with maximal Thurston-Bennequin number.
The characteristic algebra encodes the information from at least
the first- and second-order Poincar\'e-Chekanov polynomials, as we
explain in Section~\ref{ssec:polynomial}.  We will demonstrate
the efficacy of our invariant, through examples, in Section~\ref{sec:ex}.

Although the results of this section hold for links as well,
we will confine our attention to knots for simplicity, except
at the end of Section~\ref{ssec:charalgdef}.

%*********************************************************************
\subsection{Definition of the characteristic algebra}
\label{ssec:charalgdef}

The characteristic algebra can be viewed as a close relative
of the DGA homology, except that it is easier to handle in general 
than the homology itself.

\begin{definition}
Let $(A,\d)$ be a DGA over $\Z[t,t^{-1}]$, where $A = \alg{a_1,\ldots,a_n}$,
and let $I$ denote the (two-sided) ideal in $A$ generated by
$\{\d a_i\,|\,1\leq i\leq n\}$.  The {\it characteristic algebra}
$\CC(A,\d)$ is defined to be the algebra $A/I$, with grading induced 
from the grading on $A$.
\end{definition}

\begin{definition}
Two characteristic algebras $A_1/I_1$ and $A_2/I_2$ are
{\it tamely isomorphic} if we can add some number of generators to
$A_1$ and the same generators to $I_1$, and similarly for $A_2$ and $I_2$,
so that there is a tame isomorphism between $A_1$ and $A_2$ sending
$I_1$ to $I_2$.
\label{def:tame}
\end{definition}

\noindent
In particular, tamely isomorphic characteristic algebras are isomorphic
as algebras.
Strictly speaking, Definition~\ref{def:tame} only makes
sense if we interpret the characteristic algebra as a pair $(A,I)$ rather
than as $A/I$, but we will be sloppy with our notation.
Recall that we defined tame isomorphism between free algebras in
Section~\ref{ssec:invariance}.

A stabilization of $(A,\d)$, as defined in 
Section~\ref{ssec:invariance}, adds two generators $e_1,e_2$ to $A$
and one generator $e_2$ to $I$; thus $A/I$ changes by adding one
generator $e_1$ and no relations.

\begin{definition}
Two characteristic algebras $A_1/I_1$ and $A_2/I_2$
are {\it equivalent} if they are tamely isomorphic, after adding a 
(possibly different) finite number of generators (but no additional
relations) to each.
\end{definition}

\begin{theorem}
Legendrian-isotopic knots have equivalent characteristic algebras.
\label{thm:charalg}
\end{theorem}

\begin{proof}
Let $(A,\d)$ be a DGA with $A = \alg{a_1,\ldots,a_n}$.  
Consider an elementary automorphism of $A$ sending $a_j$ to $a_j + v$,
where $v$ does not involve $a_j$;
since $\d(a_j + v)$ is in $I$, it is easy to see that this automorphism
descends to a map on characteristic algebras.  We conclude that
tamely isomorphic DGAs have tamely isomorphic characteristic algebras.
On the other hand, equivalence of characteristic algebras is
defined precisely to be preserved under stabilization of DGAs.
\end{proof}

In the case of a link, we may also define the {\it characteristic
module} arising from the differential link module $(\Gamma,\d')$
introduced in Section~\ref{ssec:dgalinks}.  This is the module over
$\Z[t_1,t_1^{-1},\ldots,t_k,t_k^{-1}]$ generated by $\Gamma$,
modulo the relations
\[
v_1 (\d'a_i) v_2 = 0 \,:\, v_1 \in \Gamma_{j_1j_2}, a_i\in\Gamma_{j_2j_3},
v_2\in\Gamma_{j_3j_4}~\text{for some $j_1,j_2,j_3,j_4$}.
\]
Define equivalence of characteristic modules similarly to
equivalence of characteristic algebras, except that replacing a generator
$a_i$ by $t_{u(a_i)}^{\pm 1} a_i$ or $t_{l(a_i)}^{\pm 1} a_i$ is
allowed.  Then Legendrian-isotopic links have equivalent characteristic
modules.  An approach along these lines is used in \cite{bib:Mi}
to distinguish between particular links.
%Applications of this fact for solid-torus links are given
%in \cite{bib:NT}.

%*********************************************************************
\subsection{Relation to the Poincar\'e-Chekanov polynomial invariants}
\label{ssec:polynomial}

In this section, we work over $\Z/2$ rather than
over $\Z[t,t^{-1}]$; simply set $t=1$ and reduce modulo 2.
%(This also forces us to grade the DGA of a knot $K$ over $\Z/(2r(K))$
%rather than over $\Z$.)
Thus we consider the DGA $(A,\d)$ of a Legendrian knot $K$ over $\Z/2$,
graded over $\Z/(2r(K))$; let $\CC=A/I$ be its characteristic algebra.
%To avoid indecipherable notation, we suppress indices where possible.

We first review the definition of the Poincar\'e-Chekanov
polynomials.   The following term is taken
from \cite{bib:EFM}.

\begin{definition}
Let $(A,\d)$ be a DGA over $\Z/2$.  An algebra map $\eps\,:\,A\rightarrow
\Z/2$ is an {\it augmentation} if $\eps(1)=1$, $\eps \circ \d = 0$, and
$\eps$ vanishes for any element in $A$ of nonzero degree.
\end{definition}

Given an augmentation $\eps$ of $(A,\d)$, write $A_\eps = \ker \eps$;
%and $\d_{\eps} = \eps \circ \d \circ \eps^{-1}$;
then $\d$ maps $(A_\eps)^n$ into itself for all $n$, and thus
$\d$ descends to a map $\d^{(n)}\,:\,
A_\eps/A_\eps^{n+1} \rightarrow A_\eps/A_\eps^{n+1}$.
We can break $A_\eps/A_\eps^{n+1}$ into graded pieces
$\sum_{i\in\Z/(2r(K))} C_i^{(n)}$, where $C_i^{(n)}$ 
denotes the piece of degree $i$.
Write $\alpha_i^{(n)} = 
\dim_{\Z/2} \ker (\d^{(n)}:C_i^{(n)} \rightarrow C_{i-1}^{(n)})$ and
$\beta_i^{(n)} = \dim_{\Z/2} \im (\d^{(n)}:C_{i+1}^{(n)} \rightarrow 
C_i^{(n)})$,
so that $\alpha_i^{(n)} - \beta_i^{(n)}$ is the dimension of
the $i$-th graded piece of the homology of $\d^{(n)}$.

\begin{definition}
The {\it Poincar\'e-Chekanov polynomial of order $n$} associated
to an augmentation $\eps$ of $(A,\d)$ is
$P_{\eps,n}(\lambda) = \sum_{i\in\Z/(2r(K))} 
\left(\alpha_i^{(n)}-\beta_{i}^{(n)}\right) \lambda^i$.
\end{definition}

\noindent
Note that augmentations of a DGA do not always exist.  

The main result of this section states that we can recover
some Poincar\'e-Chekanov polynomials from the characteristic algebra.
To do this, we need one additional bit of information, besides
the characteristic algebra.

\begin{definition}
Let $\gamma_i$ be the number of
generators of degree $i$ of a DGA $(A,\d)$ graded over $\Z/(2r(K))$.  
Then the {\it degree distribution}
$\gamma\,:\,\Z/(2r(K)) \rightarrow \Z_{\geq 0}$ of $A$ is the map 
$i \mapsto \gamma_i$.
\end{definition}

\noindent
Clearly, the degree distribution can be immediately computed from a 
diagram of $K$ by calculating the degrees of the vertices of $K$.

We are now ready for the main result of this section.  Note that
the following proposition uses the {\it isomorphism} class, not
the equivalence class, of the characteristic algebra.

\begin{proposition}
The set of first- and second-order Poincar\'e-Chekanov polynomials for all 
possible augmentations
of a DGA $(A,\d)$ is determined by the isomorphism class of the 
characteristic algebra $\CC$ and the degree distribution of $A$.
\label{prop:poly}
\end{proposition}

Before we can prove Proposition~\ref{prop:poly}, we need to 
establish a few ancillary results.
Our starting point is the observation that 
there is a one-to-one correspondence between augmentations
and maximal ideals
$\langle a_1+c_1,\ldots,a_n+c_n\rangle \subset A$ containing $I$ 
and satisfying $c_i = 0$ if $\deg a_i \neq 0$.

Fix an augmentation $\eps$.  We first assume for convenience that $\eps = 0$; 
then $I \subset M$, where
$M$ is the maximal ideal $\langle a_1,\ldots,a_n \rangle$.
For each $i$, write
\[
\d a_i = \d_1 a_i + \d_2 a_i + \d_3 a_i,
\]
where $\d_1 a_i$ is linear in the $a_j$, $\d_2 a_i$ is quadratic in the
$a_j$, and $\d_3 a_i$ contains terms of third or higher order.
The following lemma writes $\d_1$ in a standard form.

\begin{lemma}
After applying a tame automorphism, we can relabel the $a_i$ as
$a_1,\ldots,a_k,b_1,\ldots,b_k$, $c_1,\ldots,c_{n-2k}$ for some $k$,
so that $\d_1 a_i = b_i$ and  $\d_1 b_i = \d_1 c_i = 0$ for all $i$.
\label{lem:poly1}
\end{lemma}

\begin{proof}
For clarity, we first relabel the $a_i$ as $\a_i$.
We may assume that the $\a_i$ are ordered so that $\d\a_i$ contains
only terms involving $\a_j$, $j < i$; see \cite{bib:Che}.
Let $i_1$ be the smallest number so that
$\d_1\a_{i_1} \neq 0$.  We can write
$\d_1\a_{i_1} = \a_{j_1} + v_1$, where $j_1 < i_1$ and the expression $v_1$
does not involve $\a_{j_1}$.  After applying the elementary isomorphism
$\a_{j_1} \mapsto \a_{j_1} + v_1$, we may assume that $v_1 = 0$
and $\d_1\a_{i_1} = \a_{j_1}$.

For any $\a_i$ such that $\d_1 \a_i$ involves $\a_{j_1}$, replace
$\a_i$ by
$\a_i + \a_{i_1}$.  Then $\d_1\a_i$ does not involve $\a_{j_1}$ unless
$i = i_1$; in addition, no $\d_1\a_i$ can involve $\a_{i_1}$, since
then $\d_1^2\a_i$ would involve $\a_{j_1}$.  Set $a_1 = \a_{i_1}$
and $b_1 = \a_{i_1}$; then $\d_1 a_1 = b_1$ and $\d_1\a_i$ does
not involve $a_1$ or $b_1$ for any other $i$.

Repeat this process with the next smallest $\a_{i_2}$ with
$\d_1\a_{i_2} \neq 0$, and so forth.  At the conclusion of this
inductive process, we obtain $a_1,\ldots,a_k,b_1,\ldots,b_k$
with $\d_1 a_i = b_i$ (and $\d_1 b_i = 0$), and the remaining
$\a_i$ satisfy $\d_1 \a_i = 0$; relabel these remaining
generators with $c$'s.
\end{proof}

Now assume that we have relabelled the generators of $A$ in 
accordance with Lemma~\ref{lem:poly1}.

\begin{lemma}
$\beta_\ell^{(1)}$ is the number of $b_j$ of degree $\ell$,
while $\beta_\ell^{(2)}-\beta_\ell^{(1)}$ is the dimension of the
degree $\ell$ subspace of the vector space generated by
\[
\{\d_2 b_i, \d_2 c_i,
a_i b_j + b_i a_j, b_i b_j, b_i c_j, c_i b_j\},
\]
where $i,j$ range over all possible indices.
\label{lem:poly2}
\end{lemma}

\begin{proof}
The statement for $\beta_\ell^{(1)}$ is obvious.  To calculate
$\beta_\ell^{(2)}-\beta_\ell^{(1)}$, note that the
image of $\d^{(2)}$ in $A/A^3$
is generated by $\d a_i = b_i + \d_2 a_i$,
$\d b_i = \d_2 b_i$, $\d c_i = \d_2 c_i$,
$\d (a_i a_j) = a_i b_j + b_i a_j$,
$\d (a_i b_j) = b_i b_j$,
$\d (b_i a_j) = b_i b_j$,
$\d (a_i c_j) = b_i c_j$,
and $\d (c_i a_j) = c_i b_j$.
\end{proof}

We wish to write $\beta_\ell^{(n)}$ in terms of $\CC$, but we first
pass through an intermediate step.
Let $N^{(n)}$ be the image of $I$ in $M/M^{n+1}$, 
and let $\delta_\ell^{(n)}$
be the dimension of the degree $\ell$ part of $N^{(n)}$.
Lemma~\ref{lem:poly4}
below relates $\beta_\ell^{(n)}$ to $\delta_\ell^{(n)}$ for
$n=1,2$.

\begin{lemma}
$\delta_\ell^{(1)}$ is the number of $b_i$ of degree $\ell$, while
$\delta_\ell^{(2)}-\delta_\ell^{(1)}$ 
is the dimension of the degree $\ell$ subspace of
the vector space generated by
\[
\{\d_2 b_i, \d_2 c_i,
a_i b_j, b_i a_j, b_i b_j, b_i c_j, c_i b_j\},
\]
where $i,j$ range over all possible indices.
\label{lem:poly3}
\end{lemma}

\begin{proof}
This follows immediately from the fact that $I$ is generated by
$\{\d a_i,\d b_i,\d c_i\}$.
\end{proof}

\begin{lemma}
$\beta_\ell^{(1)} = \delta_\ell^{(1)}$ and
$\beta_\ell^{(2)} = \delta_\ell^{(2)}
-\sum_{\ell'} \delta_{\ell'} \delta_{\ell-\ell'-1}$.
\label{lem:poly4}
\end{lemma}

\begin{proof}
We use Lemmas \ref{lem:poly2} and \ref{lem:poly3}.  The first equality
is obvious.  For the second equality, we claim that,
for fixed $i$ and $j$,
$a_i b_j$ only appears in conjunction with $b_i a_j$ in the
expressions $\d_2 b_m$ and $\d_2 c_m$, for arbitrary $m$.  
It then follows that
$\delta_\ell^{(2)} - \beta_\ell^{(2)}$ is the number of
$a_i b_j$ of degree $\ell$, which is 
$\sum_{\ell'} \delta_{\ell'} \delta_{\ell-\ell'-1}$.

To prove the claim, suppose that $\d_2 b_m$
contains a term $a_i b_j$.  Since $\d_2^2 b_m = 0$ and
$\d_2 (a_i b_j) = b_i b_j$, there must be another term in $\d_2 b_m$
which, when we apply $\d_2$, gives $b_i b_j$; but this term can only be
$b_i a_j$.  The same argument obviously holds for $\d_2 c_m$.
\end{proof}

Now let $\eps$ be any augmentation, and let $M_\eps =
\langle a_1 + \eps(a_1),\ldots,a_n + \eps(a_n)\rangle$
be the corresponding maximal ideal in $A$.
If we define $N^{(n)}$ and $\delta_\ell^{(n)}$ as above, except with
$M$ replaced by $M_\eps$, then Lemma~\ref{lem:poly4} still holds.
We are now ready to prove Proposition~\ref{prop:poly}.

\begin{proof}[Proof of Proposition \ref{prop:poly}]
Note that
\[
(M_\eps/M_\eps^{n+1})/N^{(n)} \cong (M_\eps/I)/(M_\eps/I)^{n+1}; 
\]
the characteristic algebra $\CC = A/I$ and the choice of 
augmentation $\eps$ determine the right hand side.  On the
other hand, the dimension of the degree $\ell$ part of
$M_\eps/M_\eps^{n+1}$ is $\gamma_\ell$ if $n=1$, and
$\gamma_\ell + \sum_{\ell'} \gamma_{\ell'} \gamma_{\ell-\ell'}$
if $n=2$.  It follows that we can calculate
$\{\delta_\ell^{(1)}\}$ and $\{\delta_\ell^{(2)}\}$ from $\CC$, $\eps$,
and $\gamma$.

Fix $n=1,2$.
By Lemma~\ref{lem:poly4}, we can then calculate
$\{\beta_\ell^{(n)}\}$ and hence the Poincar\'e-Chekanov polynomial
\[
P_{\eps,n}(\lambda) = \sum_\ell \left(
(\alpha_\ell^{(n)} + \beta_{\ell-1}^{(n)})
- \beta_\ell^{(n)} - \beta_{\ell-1}^{(n)} \right) \lambda^\ell.
\]
Letting $\eps$ vary over all possible augmentations yields
the proposition.
\end{proof}

The situation for higher-order Poincar\'e-Chekanov polynomials
seems more difficult; we tentatively make the following conjecture.

\begin{conjecture}
The isomorphism class of $\CC$ and
the degree distribution of $A$ determine the Poincar\'e-Chekanov 
polynomials in all orders.
\end{conjecture}

Another set of invariants, similar to the Poincar\'e-Chekanov polynomials, are
obtained by ignoring the grading of the DGA, and considering ungraded
augmentations.  In this case, the invariants are a set of
integers, rather than polynomials, in each order.  A proof
similar to the one above shows that the first- and second-order
ungraded invariants are determined by the characteristic algebra.

In practice, we apply Proposition~\ref{prop:poly} as follows.
Given two DGAs, stabilize each with the appropriate number and degrees
of stabilizations so that the two resulting DGAs have the same
degree distribution.  If these new DGAs have isomorphic
characteristic algebras, then they have the same first- and second-order
Poincar\'e-Chekanov polynomials (if augmentations exist).  If not, then we can often
see that their characteristic algebras are not equivalent, and so
the original DGAs are not equivalent.  Thus calculating characteristic
algebras often obviates the need to calculate first- and second-order
Poincar\'e-Chekanov polynomials.

Note that the {\it first-order} Poincar\'e-Chekanov polynomials 
depend only on the {\it abelianization}
of $(A,\d)$; if the procedure described above
yields two characteristic algebras whose abelianizations are
isomorphic, then the original DGAs have the same first-order
Poincar\'e-Chekanov polynomials.
On a related note, empirical evidence leads us to propose the following
conjecture, which would yield a new topological knot invariant.

\begin{conjecture}
For a Legendrian knot $K$ with maximal Thurston-Bennequin number,
the equivalence class of the abelianized characteristic algebra of
$K$, considered without grading and over $\Z$, depends only on
the topological class of $K$.
\label{conj:abelianization}
\end{conjecture}

\noindent
Here the abelianization is unsigned: $vw=wv$ for all $v,w$.

We can view the abelianization of $\CC$ in terms of algebraic
geometry.  If $\CC = (\Z/2)\langle a_1,\ldots,a_n \rangle \, / \, I$,
then the abelianization of $\CC$ gives rise to a scheme $X$ in
$\A^n$, affine $n$-space over $\Z/2$.  Theorem~\ref{thm:charalg}
immediately implies the following result.

\begin{corollary}
The scheme $X$ is a
Legendrian-isotopy invariant, up to changes of coordinates and
additions of extra coordinates (i.e., we can replace $X \subset \A^n$
by $X \times \A \subset \A^{n+1}$).
\end{corollary}

There is a conjecture about first-order Poincar\'e-Chekanov polynomials,
suggested by Chekanov, which
has a nice interpretation in our scheme picture.

\begin{conjecture}[\cite{bib:Che}]
The first-order Poincar\'e-Chekanov polynomial is independent of the
augmentation $\eps$.
\label{conj:unique}
\end{conjecture}

\noindent
Augmentations are 
simply the $(\Z/2)$-rational points in $X$, graded in the sense
that all coordinates corresponding to nonzero-degree $a_j$ are zero.
It is not hard to see that the first-order Poincar\'e-Chekanov polynomial at a 
$(\Z/2)$-rational point $p$ in $X$ is precisely
the ``graded'' codimension in $\A^n$ of $T_pX$, the tangent space
to $X$ at $p$.  The following conjecture, which we have verified
in many examples, would imply Conjecture~\ref{conj:unique}.

\begin{conjecture}
The scheme $X$ is irreducible and smooth at each $(\Z/2)$-rational point.
\end{conjecture}

\vspace{24pt}

%*********************************************************************

\section{Applications}
\label{sec:ex}

In this section, we give several illustrations of the constructions
and results from Sections~\ref{sec:frontformulation} and \ref{sec:charalg},
especially Theorems \ref{thm:polylink} and \ref{thm:charalg}.
The first three examples, all knots, both illustrate the computation of
the characteristic algebra described in Section~\ref{ssec:charalgdef},
and demonstrate its usefulness in distinguishing between Legendrian knots.
The last two examples, multi-component links, apply the techniques of 
Section~\ref{ssec:dgalinks} to conclude results about Legendrian
links.

Instead of using the full DGA over $\Z[t,t^{-1}]$
or $\Z[t_1,t_1^{-1},\ldots,t_k,t_k^{-1}]$, we will work
over $\Z/2$ by setting $t=1$ and reducing modulo 2.

%*********************************************************************
\subsection{Example 1: $\boldsymbol{6_2}$}
\label{ssec:ex1}

\begin{figure}[b]
\centering{
\includegraphics[width=0.8in,angle=270]{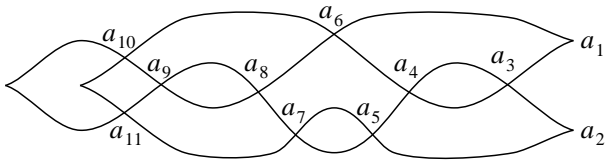}
}
\caption{The Legendrian knot $K$, of type $6_2$, with 
vertices labelled.}
\label{fig:sixtwo}
\end{figure}

Our first example answers the Legendrian mirror question of Fuchs
and Tabachnikov \cite{bib:FT}; see also \cite{bib:Ng2}.
Let the {\it Legendrian mirror} of a Legendrian knot in $\R^3$ be
the image of the knot under the involution $(x,y,z) \mapsto (x,-y,-z)$.
It is asked in \cite{bib:FT} whether a Legendrian knot with $r=0$
must always be Legendrian isotopic to its mirror.  We show that
the answer is negative by using the characteristic algebra.
Our proof is essentially identical to, but slightly cleaner than, 
the one given in \cite{bib:Ng2}; rather than using the characteristic
algebra, \cite{bib:Ng2} performs an explicit computation on the DGA
homology.

Let $K$ be the unoriented Legendrian knot given in
Figure~\ref{fig:sixtwo}, which is of knot type $6_2$, with $r=0$
and $tb=-7$.  
With vertices labelled as in Figure~\ref{fig:sixtwo}, the differential on
the DGA $(A,\d)$ for $K$ is given by
$A = \Z\langle a_1,\ldots,a_{11} \rangle$ and
\[
\begin{array}{rclcrcl}
\d a_1 &=& 1 + a_{10} a_5 a_3 & \hspace{0.5in} &
\d a_6 &=& a_{11} a_8 \\
\d a_2 &=& 1 + a_3(1+a_6a_{10}+a_{11}a_7) &&
\d a_7 &=& a_8 a_{10} \\
\d a_4 &=& a_{11} + (1 + a_6 a_{10} + a_{11} a_7) a_5 &&
\d a_9 &=& 1 + a_{10} a_{11} \\
\multicolumn{7}{c}{
\d a_5 = \d a_8 = \d a_{10} = \d a_{11} = \d a_3 = 0.}
\end{array}
\]
The ideal $I$ is generated by the above expressions;
%$I = \langle f_1,f_2,f_3,f_4,f_5,f_6 \rangle$, where
%\[
%\begin{array}{rclcrcl}
%f_1 &=& 1 + a_{10} a_5 a_3 & \hspace{0.5in} &
%f_4 &=& a_8 a_{10} \\
%f_2 &=& 1 + a_3(1+a_6a_{10}+a_{11}a_7) &&
%f_5 &=& a_{11} a_8 \\
%f_3 &=& a_{11} + (1 + a_6 a_{10} + a_{11} a_7) a_5 &&
%f_6 &=& 1 + a_{10} a_{11}.
%\end{array}
%\]
the characteristic algebra of $K$ is $\CC = A/I$.
The grading on $A$ and $\CC$ is as follows: $a_1$, $a_2$, $a_7$, 
$a_9$, and $a_{10}$ have degree 1; $a_3$, $a_4$ have degree 0; and 
$a_5$, $a_6$, $a_8$, $a_{11}$ have degree $-1$.

The characteristic algebra for the Legendrian mirror of $K$
is the same as $\CC = A/I$, but with each term in $I$ reversed.

\begin{lemma} \label{lem:charalgcomp}
We have
\[
%\begin{split}
\CC \cong (\Z/2)\langle a_1,\ldots,a_7,a_9,a_{10}\rangle \, /
%& \qquad \qquad
\langle 1+a_{10}a_5a_3,1+a_3a_{10}a_5,1+a_{10}^2a_5^2,
1+a_{10}a_5+a_6a_{10}+a_{10}a_5^2a_7 \rangle.
%\end{split}
\]
\end{lemma}

\begin{proof}
We perform a series of computations in $\CC = A/I$:
\begin{gather*}
a_8 = a_8 + (1+a_{10}a_{11}) a_8 = a_{10}(a_{11}a_8) = 0;\\
1+a_6a_{10}+a_{11}a_7 = a_{10}a_5a_3(1+a_6a_{10}+a_{11}a_7) = a_{10}a_5;\\
a_{11} = (1+a_6a_{10}+a_{11}a_7)a_5 = a_{10}a_5^2.
\end{gather*}
Substituting for $a_8$ and $a_{11}$ in the relations $\{\d a_i = 0\}$ yields the
relations in the statement of the lemma.  Conversely, given the relations
in the statement of the lemma, and setting $a_8 = 0$ and 
$a_{11} = a_{10}a_5^2$, we can recover the relations $\{\d a_i = 0\}$.
\end{proof}

Decompose $\CC$ into graded pieces $\CC = \oplus_i \CC_i$, where 
$\CC_i$ is the submodule of degree $i$.

\begin{lemma}
There do not exist $v \in \CC_{-1}, w\in \CC_1$ such that $vw = 1 \in \CC$.
\label{lem:asymmetry}
\end{lemma}

\begin{proof}
Suppose otherwise, and consider the algebra $\CC'$ obtained from 
$\CC$ by setting $a_3 = 1, a_1=a_2=a_6=a_7=a_9= 0$.  There is an obvious 
projection from $\CC$ to $\CC'$ which is an algebra map; under this
projection, $v,w$ map to $v'\in \CC'_{-1}, w'\in \CC'_1$, with
$v'w' = 1$ in $\CC'$.  But it is easy to see that
$\CC' = (\Z/2)\langle a_5,a_{10}\rangle\,/\,\langle 1+a_{10}a_5 \rangle$,
with $a_5\in \CC'_{-1}$ and $a_{10}\in \CC'_1$, and it follows that
there do not exist such $v',w'$.
\end{proof}

\begin{proposition}
$K$ is not Legendrian isotopic to its Legendrian mirror.
\end{proposition}

\begin{proof}
Let $\tilde{\CC}$ be the characteristic algebra of the Legendrian mirror
of $K$.  Since the relations in $\tilde{\CC}$ are precisely the relations
in $\CC$ reversed, Lemma~\ref{lem:asymmetry} implies that there do not
exist $v\in\tilde{\CC}_1,w\in\tilde{\CC}_{-1}$ such that
$vw=1$.  On the other hand, there certainly exist 
$v\in\CC_1,w\in\CC_{-1}$ such that $vw=1$; for instance, take
$v=a_{10}$ and $w=a_5a_3$.  Hence $\CC$ and $\tilde{\CC}$ are not isomorphic.
This argument still holds if some number of generators is added to
$\CC$ and $\tilde{\CC}$, and so $\CC$ and $\tilde{\CC}$ are not equivalent.
The result follows from Theorem~\ref{thm:charalg}.
\end{proof}

More generally, the characteristic algebra technique seems to be an effective
way to distinguish between some knots and their Legendrian mirrors;
cf.\ Section~\ref{ssec:ex2}.
%see Remark~\ref{rmk:sevenfour} for another example.
Note that Poincar\'e-Chekanov polynomials of any order can 
never tell between a knot
and its mirror, since, as mentioned above, the differential for a mirror
is the differential for the knot, with each monomial reversed.

%*********************************************************************
\subsection{Example 2: $\boldsymbol{7_4}$}
\label{ssec:ex2}

\begin{figure}
\centering{
\includegraphics[width=2.5in,angle=270]{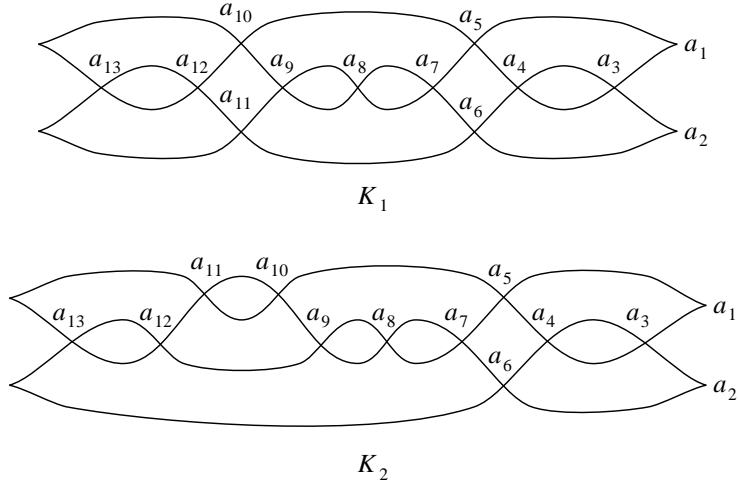}
}
\caption{
The Legendrian knots $K_1$ and $K_2$, of type $7_4$,
with vertices labelled.
}
\label{fig:sevenfour}
\end{figure}

Our second example shows that the characteristic algebra is effective
even when Poincar\'e-Chekanov polynomials do not exist.  In addition,
this section and the next provide the first examples, known to the author,
in which the DGA grading is not needed to distinguish between knots.

Consider the Legendrian knots $K_1$, $K_2$ shown in 
Figure~\ref{fig:sevenfour}; both are of smooth type $7_4$, with
$r=0$ and $tb=1$.  We will show that $K_1$ and $K_2$ are not Legendrian
isotopic.
%  We present this example before the $6_3$ and $7_2$ examples
%of Sections \ref{sec:ex3new} because the
%algebra is a bit simpler in this case.

The differential on the DGA for $K_1$ is given by
\begin{alignat*}{2}
\d a_1 &= 1 + a_8a_{13}a_3 & \qquad \d a_6 &= (1+a_7a_8) a_{13} \\
\d a_2 &= 1 + a_3a_{13}a_8 & \d a_9 &= a_{13}a_{11} + a_{10}a_{13} \\
\d a_4 &= a_5a_8a_{13} + a_{13}a_8a_6 & \d a_{10} &= 1 + a_{13}a_{12} \\
\d a_5 &= a_{13} (1+a_8a_7) & \d a_{11} &= 1 + a_{12}a_{13};
\end{alignat*}
the differential for $K_2$ is given by
\begin{align*}
\d a_1 &= 1 + (1+a_8a_9+a_8a_{13}+a_{12}a_{13}+a_8a_9a_{12}a_{13}+a_8a_{10}a_{11}a_{13}) a_3 \\
\d a_2 &= 1 + a_3a_{13}a_8 \\
\d a_4 &= a_{13}a_8a_6 + a_{11}a_{13} 
+ a_5 (1+a_8a_9+a_8a_{13}+a_{12}a_{13}+a_8a_9a_{12}a_{13}+a_8a_{10}a_{11}a_{13}) \\
\d a_5 &= a_{13} (1+a_8a_7) \\
\d a_6 &= a_7+a_7a_{12}a_{13} + 
(1+a_7a_8)(a_9+a_{13}+a_9a_{12}a_{13}+a_{10}a_{11}a_{13}) \\
\d a_9 &= a_{10}a_{13} \\
\d a_{11} &= 1 + a_{13}a_{12}.
\end{align*}

Denote the characteristic algebras of $K_1$ and $K_2$ by $\CC_1=A/I_1$ and
$\CC_2=A/I_2$, respectively; here 
$A = (\Z/2) \langle a_1,\ldots,a_{13} \rangle$,
and $I_1$ and $I_2$ are generated by the respective expressions above.

\begin{lemma}
We have
\[
\CC_1 \cong (\Z/2) \langle a_1,\ldots,a_5,a_8,a_9,a_{10},a_{13} \rangle
\, / \, \langle 1 + a_3a_8a_{13}, a_3a_8 + a_8a_3, a_3a_{13}+a_{13}a_3,
a_8a_{13}+a_{13}a_8 \rangle .
\]
\end{lemma}

\begin{proof}
Similar to the proof of Lemma~\ref{lem:charalgcomp}.
\end{proof}

\begin{lemma} \label{lem:sevenfour1}
There is no expression in $\CC_1$ which is invertible from one side
but not from the other.
\end{lemma}

\begin{proof}
It is clear that the only expressions in $\CC_1$ which are invertible
from either side are
products of some number of $a_3$, $a_8$, and $a_{13}$, with inverses
of the same form.  Since $a_3,a_8,a_{13}$ all commute, the lemma
follows.
\end{proof}

\begin{lemma}
In $\CC_2$, $a_{13}$ is invertible from the right but not from the left.
\label{lem:sevenfour2}
\end{lemma}

\begin{proof}
Since $a_{13}a_{12}=1$, $a_{13}$ is certainly invertible from the right.
Now consider adding to $\CC_2$ the relations
$a_3=1$, $a_7=a_{13}$, $a_8=a_{12}$, and
$a_i = 0$ for all $i$ not previously mentioned.  
The resulting algebra is isomorphic to 
$(\Z/2)\langle a_{12},a_{13} \rangle \, / \, \langle 1 + a_{13}a_{12} \rangle$,
in which $a_{13}$ is not invertible from the left.  We conclude that
$a_{13}$ is not invertible from the left in $\CC_2$ either, as desired.
\end{proof}

\begin{proposition}
The Legendrian knots $K_1$ and $K_2$ are not Legendrian isotopic.
\end{proposition}

\begin{proof}
From Lemmas \ref{lem:sevenfour1} and \ref{lem:sevenfour2}, 
$\CC_1$ and $\CC_2$ are not equivalent.
\end{proof}

Although $\CC_1$ and $\CC_2$ are not equivalent, one may compute
that their abelianizations are isomorphic; 
cf.\ Conjecture~\ref{conj:abelianization}.  It is also easy to check that
$K_1$ and $K_2$ have no augmentations, and hence no
Poincar\'e-Chekanov polynomials.

The computation from the proof of Lemma~\ref{lem:sevenfour2} also
demonstrates that $K_2$ is not Legendrian isotopic to its Legendrian mirror;
we may use the same argument as in Section~\ref{ssec:ex1}, along
with the fact that $a_{12}$ and $a_{13}$ have degrees $2$ and $-2$,
respectively, in $\CC_2$.  By contrast, we see from inspection that 
$K_1$ is the same as its Legendrian mirror.

%*********************************************************************
\subsection{Example 3: $\boldsymbol{6_3}$ and $\boldsymbol{7_2}$}
\label{ssec:ex3}

\begin{figure}
\centering{
\includegraphics[width=2.5in,angle=270]{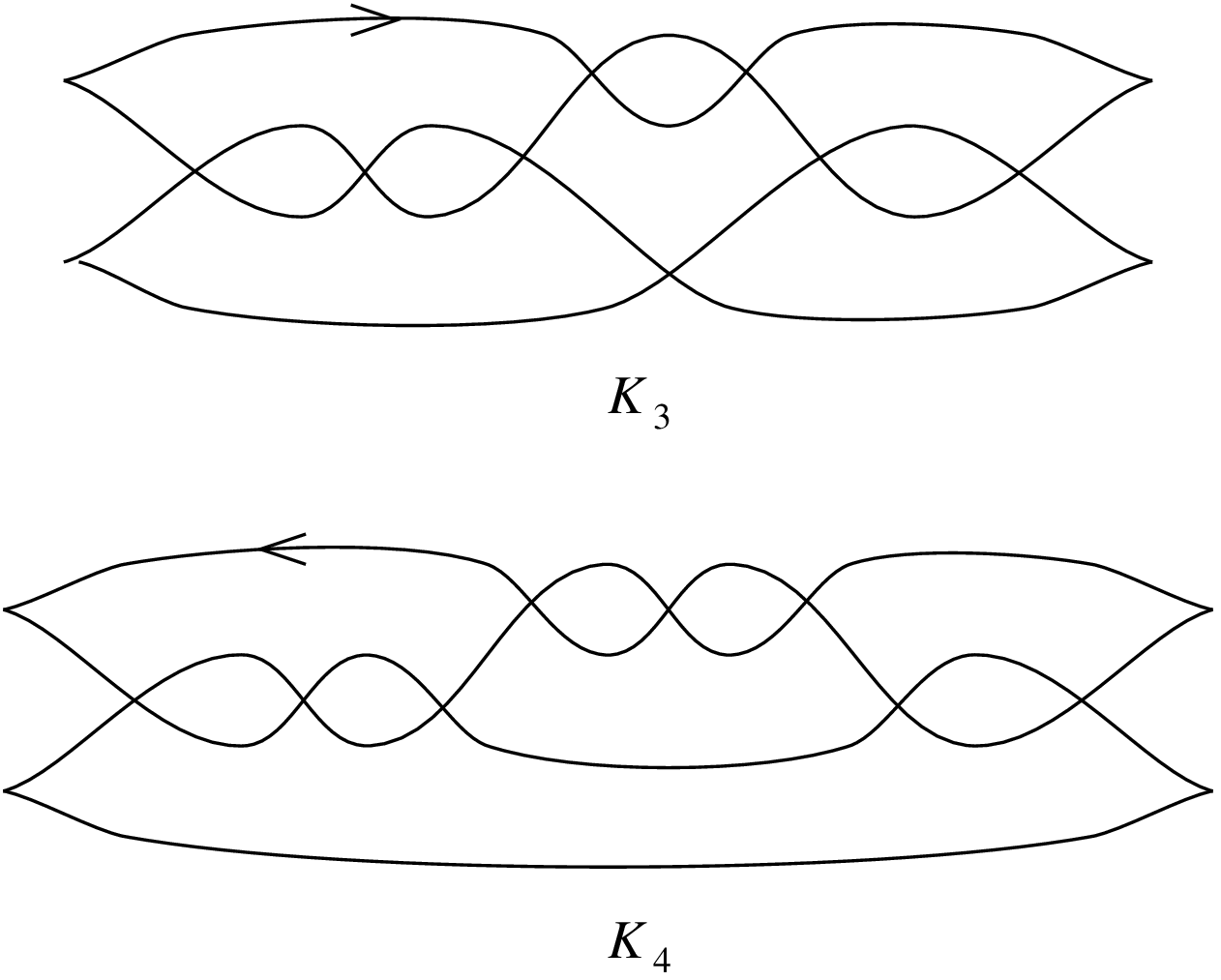}

\vspace{0.2in}

\includegraphics[width=2.5in,angle=270]{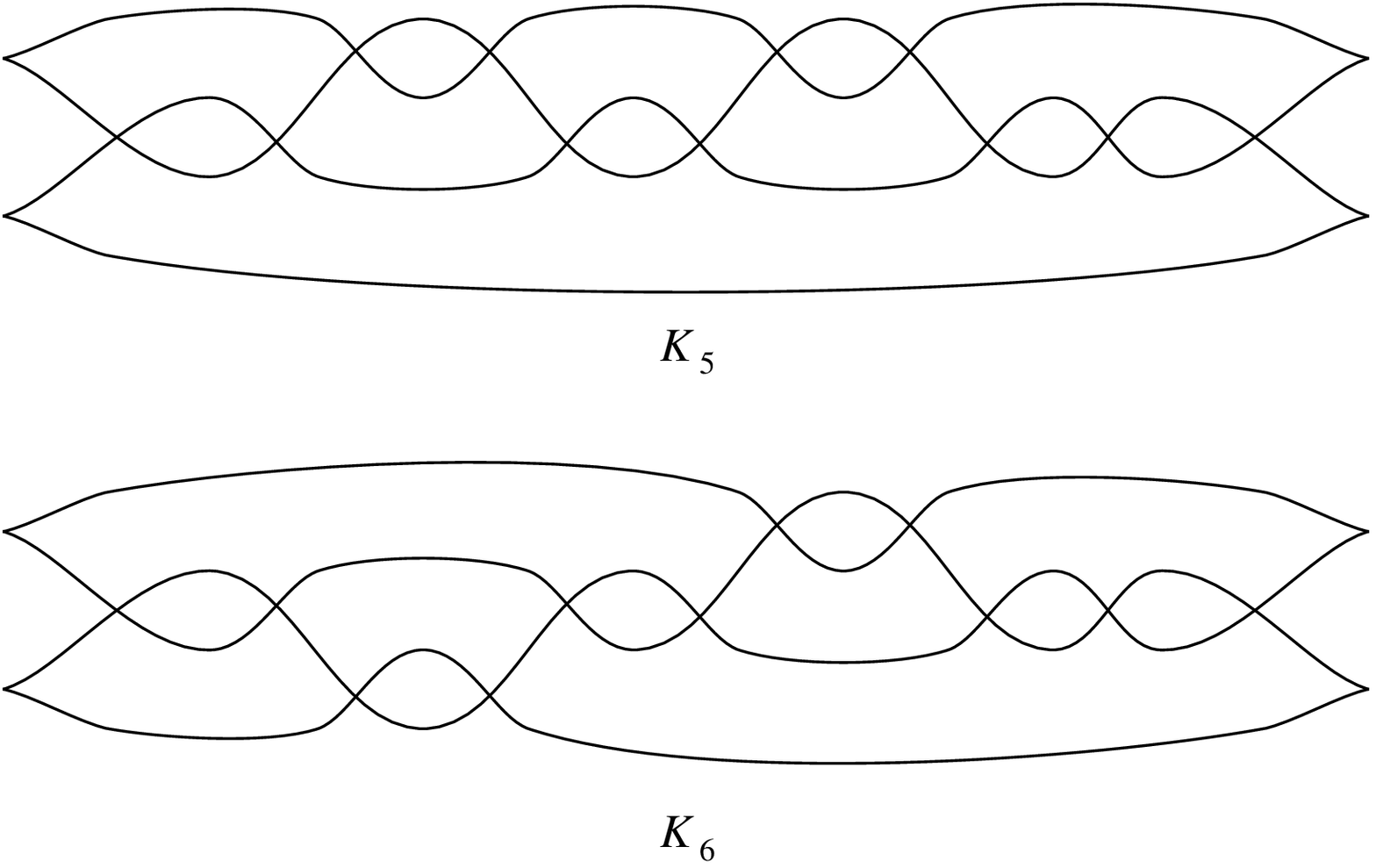}
}
\caption{
The oriented Legendrian knots $K_3$ and $K_4$, of type $6_3$,
and the unoriented knots $K_5$ and $K_6$, of type $7_2$.
}
\label{fig:sixthree}
\end{figure}

In a manner entirely analogous to Section~\ref{ssec:ex2}, we can
prove that many other pairs of Legendrian knots are not Legendrian
isotopic.  For example, consider the knots in Figure~\ref{fig:sixthree}:
$K_3$ and $K_4$, of smooth type $6_3$, with $r=1$ and $tb=-4$, 
and $K_5$ and $K_6$, of smooth type $7_2$, with $r=0$ and $tb=1$.

\begin{proposition} \label{prop:combinedex}
$K_3$ and $K_4$ are not Legendrian isotopic; $K_5$ and $K_6$ are
not Legendrian isotopic.
\end{proposition}

The proof of Proposition~\ref{prop:combinedex}, which involves
computations on the characteristic algebra along the lines of
Section~\ref{ssec:ex2}, is omitted here, but can be found in
\cite{bib:thesis}.
The $6_3$ examples are the first, known to the author, of two knots
with nonzero rotation number which have the same classical invariants
but are not Legendrian isotopic.  The first-order 
Poincar\'e-Chekanov polynomial
fails to distinguish between either the $6_3$ or the $7_2$ knots;
$K_3$ and $K_4$ have no augmentations, while $K_5$ and $K_6$ both
have first-order polynomial $\lambda + 2$.

%*********************************************************************
\subsection{Example 4: triple of the unknot}
\label{ssec:ex4}

In this section, we rederive a result of \cite{bib:Mi} by using the link 
grading from Section~\ref{ssec:dgalinks}.  Our proof is different from 
the ones in \cite{bib:Mi}.

\begin{definition}[\cite{bib:Mi}]
Given a Legendrian knot $K$, let the {\it $n$-copy} of
$K$ be the link consisting of $K$, along with $n-1$ copies of $K$
slightly perturbed in the transversal direction.  In the front projection,
the $n$-copy is simply $n$ copies of the front of $K$, differing from 
each other by small shifts in the $z$ direction.  We will call the 
$2$-copy and $3$-copy the {\it double} and {\it triple}, respectively.
\label{def:ncopy}
\end{definition}

Let $L=(L_1,L_2,L_3)$ be the unoriented triple of the usual ``flying-saucer''
unknot; this is the unoriented version of the link shown in 
Figure~\ref{fig:triple}.

\begin{proposition}[\cite{bib:Mi}]
The unoriented links $(L_1,L_2,L_3)$ and $(L_2,L_1,L_3)$ are not
Legendrian isotopic.
\end{proposition}

\begin{proof}
In Example~\ref{ex:polytriple}, we have already calculated the
first-order Poincar\'e-Chekanov polynomials for $(L_1,L_2,L_3)$, once
we allow the grading $(\rho_1,\rho_2)$ to range in
$(\frac{1}{2}\Z)^2$.
The polynomials for the link $(L_2,L_1,L_3)$ and grading
$(\sigma_1,\sigma_2) \in (\frac{1}{2}\Z)^2$ are identical,
except with the indices $1$ and $2$ reversed.
%\begin{alignat*}{3}
%P^{\eps,1}_{11}(\lambda) &= \lambda & \qquad 
%P^{\eps,1}_{21}(\lambda) &= \lambda^{-1 + 2\sigma_1 - 2\sigma_2} & \qquad
%P^{\eps,1}_{31}(\lambda) &= \lambda^{1 - 2\sigma_2} \\
%P^{\eps,1}_{12}(\lambda) &= \lambda^{1 + 2\sigma_2 - 2\sigma_1} & \qquad 
%P^{\eps,1}_{22}(\lambda) &= \lambda & \qquad
%P^{\eps,1}_{32}(\lambda) &= \lambda^{1 - 2\sigma_1} \\
%P^{\eps,1}_{13}(\lambda) &= \lambda^{-1 + 2\sigma_2} & \qquad 
%P^{\eps,1}_{23}(\lambda) &= \lambda^{-1 + 2\sigma_1} & \qquad
%P^{\eps,1}_{33}(\lambda) &= \lambda.
%\end{alignat*}
It is easy to compute that there is no choice of 
$\rho_1,\rho_2,\sigma_1,\sigma_2$ for which these polynomials
coincide with the polynomials for $(L_1,L_2,L_3)$ given in 
Example~\ref{ex:polytriple}.  The result now follows from
Theorem~\ref{thm:polylink}.
\end{proof}

%*********************************************************************
\subsection{Example 5: other links}
\label{ssec:ex5}

\begin{figure}

\vspace{0.2in}

\centering{
\includegraphics[width=2in,angle=270]{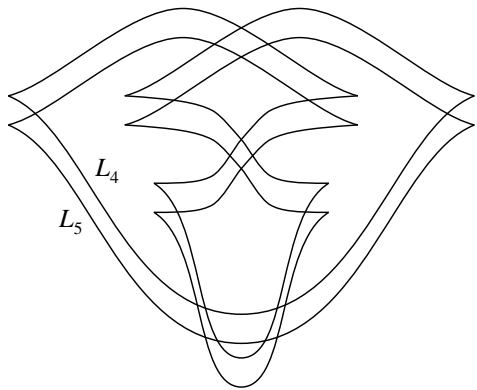}
}
\caption{
The double of the figure eight knot from Figure~\ref{fig:figure8}.
}
\label{fig:figure8double}
\end{figure}

In this section, we give two examples of other links which can be distinguished
using our techniques.  The proofs, which are simple
and can be found in \cite{bib:thesis}, use the
Poincar\'e-Chekanov polynomial and Theorem~\ref{thm:polylink}, 
as in Section~\ref{ssec:ex4}.

Let $(L_4,L_5)$ be the unoriented double of the figure
eight knot, shown in Figure~\ref{fig:figure8double}, and let $(L_5,L_4)$
be the same link, but with components interchanged.  The following
result answers a question from \cite{bib:Mi} about whether
there is an unoriented knot whose double is not isotopic
to itself with components interchanged.

\begin{proposition} \label{prop:figure8double}
The unoriented links $(L_4,L_5)$ and $(L_5,L_4)$ are not Legendrian
isotopic.
\end{proposition}

\begin{figure}
\centering{
\includegraphics[width=1.5in,angle=270]{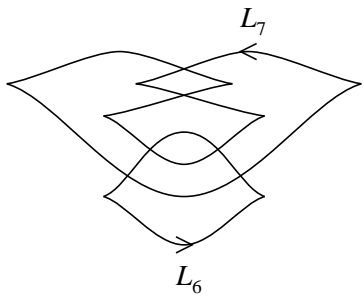}
}
\caption{
The oriented Whitehead link.
}
\label{fig:whitehead}
\end{figure}

Our other example, in which orientation is important, 
is the oriented Whitehead link $(L_6,L_7)$ shown in 
Figure~\ref{fig:whitehead}.  Let $-L_j$ denote $L_j$ with reversed 
orientation.  
By playing with the diagrams, one can
show that $(L_6,L_7)$, $(L_7,-L_6)$, $(-L_6,-L_7)$, and
$(-L_7,L_6)$ are Legendrian isotopic, as are
$(-L_6,L_7)$, $(-L_7,-L_6)$, $(L_6,-L_7)$, and $(-L_7,-L_6)$.
It is also the case that these two families are smoothly isotopic
to each other.  By contrast, we have the following result.

\begin{proposition}
The oriented links $(L_6,L_7)$ and $(-L_7,L_6)$ are not 
Legendrian isotopic.
\label{prop:whiteheadlink}
\end{proposition}

\vspace{24pt}

%*********************************************************************
\section*{Acknowledgments}

The author is grateful to Tom Mrowka, John Etnyre, Kiran
Kedlaya, Josh Sabloff, and Lisa
Traynor for useful discussions, to Isadore Singer for his
encouragement and support, and to the American Institute of Mathematics
for sponsoring the Low-Dimensional Contact Geometry program in
the fall of 2000, which was invaluable for this paper.

\vspace{24pt}

%*********************************************************************

\vspace{24pt}

\end{document}